\documentclass[aos,MSNbibl,dvips]{arximspdf}

\doi{10.1214/13-AOS1150} 
\volume{41}
\issue{5}
\pubyear{2013}
\firstpage{2359}
\lastpage{2390}

\makeatletter
\newcommand{\rrvert}{\vert}
\newcommand{\llvert}{\vert}

\newcommand{\iint}{\int\!\!\int}

\newtheorem{theorem}{Theorem}
\newtheorem{lemma}{Lemma}

\newproclaim{assumption}{Assumption}

\newtheorem{proposition}{Proposition}

\newproclaim{remark}{Remark}

\newcommand{\spn}{\operatorname{span}}
\newcommand{\op}{\mathrm{op}}
\newcommand{\TV}{\mathrm{TV}}
\newcommand{\inte}{\mathrm{int}}
\newcommand{\Var}{\operatorname{Var}}
\newcommand{\tr}{\operatorname{tr}}

\newcommand{\E}{\mathbb{E}}
\newcommand{\Prob}{\mathbb{P}}
\newcommand{\En}{\mathbb{E}_{n}}
\newcommand{\R}{\mathbb{R}}
\newcommand{\C}{\mathcal{C}}
\newcommand{\D}{\mathcal{D}}
\newcommand{\G}{\mathcal{G}}

\makeatother

\begin{document}
\begin{frontmatter}

\title{Quasi-Bayesian analysis of nonparametric instrumental variables models\thanksref{T1}}
\runtitle{Quasi-Bayes for NPIV}

\thankstext{T1}{Supported by the Grant-in-Aid for Young Scientists (B) (25780152) from the JSPS.}

\begin{aug}
\author[A]{\fnms{Kengo} \snm{Kato}\corref{}\ead[label=e1]{kkato@e.u-tokyo.ac.jp}}
\runauthor{K. Kato}
\affiliation{University of Tokyo}
\address[A]{Graduate School of Economics\\
University of Tokyo \\
7-3-1 Hongo, Bunkyo-ku\\
Tokyo 113-0033\\
Japan\\
\printead{e1}} 
\end{aug}

\received{\smonth{12} \syear{2012}}
\revised{\smonth{7} \syear{2013}}

%
\begin{abstract}
This paper aims at developing a quasi-Bayesian analysis of the
nonparametric instrumental variables model, with a focus on the
asymptotic properties of quasi-posterior distributions. In this paper,
instead of assuming a \mbox{distributional} assumption on the data generating
process, we consider a quasi-likelihood induced from the conditional
moment restriction, and put priors on the function-valued parameter. We
call the resulting posterior quasi-posterior, which corresponds to
``Gibbs posterior'' in the literature. Here we focus on priors
constructed on slowly growing finite-dimensional sieves. We derive
rates of contraction and a nonparametric Bernstein--von Mises type
result for the quasi-posterior distribution, and rates of convergence
for the quasi-Bayes estimator defined by the posterior expectation. We
show that, with priors suitably chosen, the quasi-posterior
distribution (the quasi-Bayes estimator) attains the minimax optimal
rate of contraction (convergence, resp.). These results greatly
sharpen the previous related work.
\end{abstract}

%

\end{frontmatter}

\section{Introduction}\label{sec1}
\label{secintro}

\subsection{Overview}\label{sec1.1}

Let $(Y,X,W)$ be a triplet of scalar random variables, where $Y$ is a
dependent variable, $X$ is an endogenous variable and $W$ is an
instrumental variable. Without loosing much generality, we assume that
the support of $(X,W)$ is contained in $[0,1]^{2}$.
The support of $Y$ may be unbounded.
We consider the nonparametric instrumental variables (NPIV) model of
the form
%
\begin{equation}\label{npiv}
\E[ Y \mid W ] = \E\bigl[ g_{0}(X) \mid W \bigr],
\end{equation}
where $g_{0}\dvtx  [0,1] \to\R$ is an unknown structural function of interest.
Alternatively, we can write the model in a more conventional form
\[
Y = g_{0}(X) + U, \E[ U \mid W ] = 0,
\]
where $X$ is potentially correlated with $U$ and hence $\E[ U \mid X ]
\neq0$.

A model of the form (\ref{npiv}) is of principal importance in
econometrics (see \mbox{\cite{HH05,H11b}}). From a statistical perspective,
the problem of recovering the structural function $g_{0}$ is
challenging since it is an \emph{ill-posed} inverse problem with an
additional difficulty of \emph{unknown} operator [$K$ in (\ref{npiv2})
ahead]. Statistical inverse problems, including the current problem,
have attracted considerable interests in statistics and econometrics
(see, e.g., \cite{CFR07,C08}). For mathematical background of inverse
problems, we refer to \cite{K99}.

To see that the problem of recovering the structural function $g_{0}$
is an ill-posed inverse problem,
suppose that $(X,W)$ has a square-integrable joint density
$f_{X,W}(x,w)$ on $[0,1]^{2}$ and denote by $f_{W}(w)$ the marginal
density of $W$. Define the linear operator $K\dvtx L_{2}[0,1] \to
L_{2}[0,1]$ by
\[
(Kg) (w) = \E\bigl[ g(X) \mid W=w \bigr] f_{W}(w) = \int g(x)
f_{X,W}(x,w) \,dx.
\]
Then the NPIV model (\ref{npiv}) is equivalent to the operator equation
%
\begin{equation}\label{npiv2}
Kg_{0} = h,
\end{equation}
where $h(w) = \E[ Y \mid W=w ] f_{W}(w)$. Suppose that $K$ is injective
to guarantee identification of $g_{0}$.\setcounter{footnote}{1}\footnote{This global
identification condition is, however, not a trivial assumption; see the
discussion after Assumption \ref{A2} in Section \ref{secassumption} as
well as the last paragraph in the next subsection.} The problem is
that, even though $K$ is injective, its inverse $K^{-1}$ is not
$L_{2}$-continuous since $K$ is Hilbert--Schmidt (as $f_{X,W}(x,w)$ is
square integrable on $[0,1]^{2}$) and hence the $l$th largest singular
value, denoted by $\kappa_{l}$, is approaching zero as $l \to\infty$
(see, e.g., \cite{Y80}). In this sense, the problem of recovering
$g_{0}$ from $h$ is ill-posed.

Approaches to estimating the structural function $g_{0}$ are roughly
classified into two types: the method involving the Tikhonov
regularization \cite{HH05,DFFR10} and the sieve-based method
\cite{NP03,AC03,BCK07,H11a}.\footnote{The sieve-method is further
classified into two types: the method using slowly growing
finite-dimensional sieves with no or light penalties where the
dimensions of sieves play the role of regularization, and the method
using large-dimensional sieves with heavy penalties where the penalty
terms play the role of regularization (see \cite{CP11}).} The minimax
optimal rates of convergence in estimating $g_{0}$ are established in
\cite{HH05,CR11}, and they are achieved by the estimators proposed in
\cite{HH05,BCK07} under their respective assumptions. All the above
mentioned studies are, however, from a purely frequentist perspective.
Little is known about the theoretical properties of Bayes or
quasi-Bayes analysis of the NPIV model. Exceptions are
\cite{FS10a,FS10b,FS13,LJ11}.

This paper aims at developing a quasi-Bayesian analysis of the NPIV
model, with a focus on the asymptotic properties of quasi-posterior
distributions. The approach taken is quasi-Bayes in the sense that it
neither needs to assume any specific distribution of $(Y,X,W)$, nor has
to put a nonparametric prior on the unknown likelihood function. The
analysis is then based upon a quasi-likelihood induced from the
conditional moment restriction. The quasi-likelihood is constructed by
first estimating the conditional moment function $m(\cdot,g) = \E[ Y
- g(X) \mid W =\cdot]$ in a nonparametric way, and taking $\exp\{
-(1/2) \sum_{i=1}^{n} \hat{m}^{2}(W_{i},g)\}$ as if it were a
likelihood of $g$.
For this quasi-likelihood, we put a prior on the function-valued
parameter $g$. By doing so, formally, the posterior distribution for
$g$ may be defined, which we call ``quasi-posterior distribution.''
This posterior corresponds to what \cite{JT08} called ``Gibbs
posterior,'' and has a substantial interpretation (see Proposition \ref
{prop0} ahead). The quasi-Bayesian approach in this paper builds upon
\cite{CH03} where the dimension of the parameter of interest is finite
and fixed.

We focus here on priors constructed on slowly growing
finite-dimensional sieves (called ``sieve or series priors''), where
the dimensions of the sieve spaces (which grow with the sample size)
play the role of regularization to deal with the problem of
ill-posedness. Potentially, there are several choices in sieve spaces,
but we choose to use wavelet bases to form sieve spaces. Wavelet bases
are useful to treat smoothness function classes such as
H\"{o}lder--Zygmund and Sobolev spaces in a unified and convenient way.
We also use wavelet series estimation of the conditional moment
function.\footnote{This does not rule out the use of other bases such
as the Fourier and Hermite polynomial bases. See Remark \ref{bases}.}

Under this setup, we study the asymptotic properties of the
quasi-posterior distribution.
The results obtained are summarized as follows. First, we derive rates
of contraction for the quasi-posterior distribution and establish
conditions on priors under which the minimax optimal rate of
contraction is attained.
Here the contraction is stated in the standard $L_{2}$-norm.
Second, we show asymptotic normality of the quasi-posterior of the
first $k_{n}$ generalized Fourier coefficients, where $k_{n} \to\infty
$ is the dimension of the sieve space.
This may be viewed as a nonparametric Bernstein--von Mises type result
(see \cite{V98}, Chapter 10, for the classical Bernstein--von Mises
theorem for regular parametric models).
Third, we derive rates of convergence of the quasi-Bayes estimator
defined by the posterior expectation and show that under some
conditions it attains the minimax optimal rate of convergence.
Finally, we give some specific sieve priors for which the
quasi-posterior distribution (the quasi-Bayes estimator) attains the
minimax optimal rate of contraction (convergence, resp.).
These results greatly sharpen the previous work of, for example, \cite
{LJ11}, as we will review below.

\subsection{Literature review and contributions}\label{sec1.2}

Closely related are \cite{FS10b} and \cite{LJ11}. The former paper
worked on the reduced form equation $Y = \E[ g_{0}(X) \mid W ] + V$
with $V = U+g_{0}(X)-\E[ g_{0}(X) \mid W ]$ and assumed $V$ to be
normally distributed. They considered a Gaussian prior on $g$, and the
posterior distribution is also Gaussian (conditionally on the variance
of $V$). They proposed to ``regularize'' the posterior and studied the
asymptotic properties of the ``regularized'' posterior distribution and
its expectation. Clearly, the present paper largely differs from \cite
{FS10b} in that (i) we do not assume normality of the ``error''; (ii)
roughly speaking, Florens and Simoni's method is tied with the Tikhonov
regularization method, while ours is tied with the sieve-based method
with slowly growing sieves.
We note the settings of \cite{FS10a,FS13} are largely different from
the present paper; moreover in the NPIV example, some high-level
conditions on estimated operators are assumed in \mbox{\cite{FS10a,FS13}},
and hence they are not directly comparable to the present paper. Liao
and Jiang
\cite{LJ11} developed an important unified framework in estimating
conditional moment restriction models based on a quasi-Bayesian
approach, and their scope is more general than ours.
They analyzed NPIV models in detail in their Section~4. Their posterior
construction is similar to ours such as the use of sieve priors, but
differs from ours in detail.
For example, \cite{LJ11} transformed the conditional moment
restriction into unconditional moment restrictions with increasing
number of restrictions. On the other hand, we directly work on the
conditional moment restriction, although whether Liao and Jiang's
approach will lose any efficiency in the frequentist sense is not
formally clear.

Importantly and substantially, neither \cite{FS10b} nor \cite{LJ11}
established sharp contraction rates for their (quasi-)posterior
distributions, nor asymptotic normality results.
It is unclear whether Florens and Simoni's \cite{FS10b} rates (in
their Theorem 2) are optimal, since their assumptions are substantially
different from the past literature such as \cite{HH05} and \cite
{CR11}; moreover, strictly speaking \cite{FS10b} did not formally
derive contraction rates for their regularized posterior when the
operator is unknown (note that \cite{FS10a,FS13}, though not directly
comparable to the present paper, also did not formally derive posterior
contraction rates in the NPIV example). Liao and Jiang \cite{LJ11}
only established posterior consistency.
Here we focus on a simple but important model, and establish the
sharper asymptotic results for the quasi-posterior distribution.
Notably, a wide class of (finite dimensional) sieve priors is shown to
lead to the optimal contraction rate.
Moreover, in \cite{LJ11}, a point estimator of the structural function
is not formally analyzed.
Hence, the primal contribution of this paper is to considerably deepen
the understanding of the asymptotic properties of the quasi-Bayesian
procedure for the NPIV model.

The present paper deals with a quasi-Bayesian analysis of an
infinite-dimen\-sional model. The literature on theoretical studies of
Bayesian analysis of infinite-dimensional models is large. See
\cite{GGV00,SW01,GR03,KV06,GV07,GN11} for general contraction rates
\mbox{results} for posterior distributions in infinite-dimensional models.
Note that these results do not directly apply to our case: the proof of
the main general theorem (Theorem \ref{thm1}) depends on the
construction of suitable ``tests'' (see the proof of Proposition
\ref{propA1}), but how to construct such tests in a specific problem in
a nonlikelihood framework is not trivial, especially in the current
NPIV model where we have to deal with the ill-posedness of inverse
problem. Moreover, Proposition \ref{propA1} alone is not sufficient for
obtaining sharp contraction rates and an additional work is needed (see
the proof of Theorem \ref{thm1}).

There is also a large literature on the Bayesian analysis of
(ill-posed) inverse problems.
One stream of research on this topic lies in the applied mathematics
literature; see \cite{S10} and references therein.
However, their models and scopes are substantially different from those
of the present paper; for example, \cite{HP07,HP09} considered
(ill-conditioned) finite-dimensional linear regression models with
Gaussian errors and priors, and contractions rates of posterior
distributions are not formally studied there.
In the statistics literature, we may refer to \cite
{C93,KVV11,KVV13,ALS13,KSVV12} (in addition to \cite
{LJ11,FS10a,FS10b,FS13} that are already discussed), although their
results are not applicable to the analysis of NPIV models because of
its particular structure (i.e., especially the operator $K$ is unknown,
and non-Gaussian ``errors'' and priors are allowed).
Hence the present paper provides a further contribution to the Bayesian
analysis of ill-posed inverse problems.

Our asymptotic normality result builds upon the previous work on
asymptotic normality of (quasi-)posterior distributions for models with
increasing number of parameters \cite
{G99,G00,BC09a,BC09b,BG09,CG10,B11}. Related is \cite{B11}, in which
the author established Bernstein--von Mises theorems for Gaussian
regression models with increasing number of regressors and improved
upon the earlier work of \cite{G99} in several aspects. Reference \cite{B11}
covered nonparametric models by taking into account modeling bias in
the analysis. However, none of these papers covered the NPIV model, nor
more generally linear inverse problems.

Finally, while we here assume injectivity of the operator $K$ in (\ref
{npiv2}), as one of anonymous referees pointed out, this condition is
not a trivial assumption (see also the discussion after Assumption \ref
{A2} in Section \ref{secassumption}), and there are a number of works
that relax the injectivity assumption and explore partial
identification approach, such as \cite{S12,LJ11,KY11} and \cite
{CP11}, Appendix A.

\subsection{Organization and notation}\label{sec1.3}

The remainder of the paper is organized as follows.
Section \ref{sec2} gives an informal discussion of the quasi-Bayesian analysis
of the NPIV model.
Section \ref{sec3} contains the main results of the paper where general theorems
on contraction rates and asymptotic normality for quasi-posterior
distributions, as well as convergence rates for quasi-Bayes estimators,
are stated.
Section \ref{sec4} analyzes some specific sieve priors.
Section~\ref{sec5} contains\vadjust{\goodbreak} the proofs of the main results.
Section \ref{sec6} concludes with some further discussions.
The Appendix contains some omitted technical results.
Because of the space limitation, the Appendix is contained in the
supplemental file~\cite{K13}.

\textit{Notation}:
For any given (random or nonrandom, scalar or vector) sequence $\{
z_{i} \}_{i=1}^{n}$, we use the notation $\En[ z_{i} ]=n^{-1}\sum
_{i=1}^{n} z_{i}$,
which should be distinguished from the population expectation $\E
[\cdot]$. For any vector $z$, let $z^{\otimes2} = z z^{T}$ where
$z^{T}$ is the transpose of $z$. For any two sequences of positive
constants $r_{n}$ and $s_{n}$, we write
$r_{n} \lesssim s_{n}$ if the ratio $r_{n}/s_{n}$ is bounded, and
$r_{n} \sim s_{n}$ if $r_{n} \lesssim s_{n}$ and $s_{n} \lesssim
r_{n}$. Let $L_{2}[0,1]$ denote the usual $L_{2}$ space with respect to
the Lebesgue measure for functions defined on $[0,1]$. Let $\| \cdot\|
$ denote the $L_{2}$-norm, that is, $\| f \|^{2} = \int_{0}^{1}
f^{2}(x) \,dx$. The inner product in $L_{2}[0,1]$ is denoted by $\langle
\cdot, \cdot\rangle$, that is, $\langle f,g \rangle= \int_{0}^{1}
f(x) g(x) \,dx$. Let $C[0,1]$ denote the metric space of all continuous
functions on $[0,1]$, equipped with the uniform metric. The Euclidean
norm is denoted by $\| \cdot\|_{\ell^{2}}$. For any matrix $A$, let
$s_{\min}(A)$ and $s_{\max}(A)$ denote the minimum and maximum
singular values of $A$, respectively. Let $\| A \|_{\op}$ denote the
operator norm of a matrix $A$ [i.e., $\| A \|_{\op} = s_{\max} (A)$].
Denote by $dN(\mu,\Sigma)(x)$ the density of the multivariate normal
distribution with mean vector $\mu$ and covariance matrix $\Sigma$.

\section{Quasi-Bayesian analysis: Informal discussion}\label{sec2}

In this section, we outline a quasi-Bayesian analysis of the NPIV model
(\ref{npiv}). The discussion here is informal. The formal discussion
is given in Section \ref{secmain}.

Let $\G$ be a parameter space (say, some smoothness class of
functions, such as a H\"{o}lder--Zygmund or Sobolev space), for which we
assume $g_{0} \in\G$. We assume that $\mathcal{G}$ is at least
contained in $C[0,1]$: $\mathcal{G} \subset C[0,1]$. Define the
conditional moment function as $m(W,g) = \E[ Y - g(X) \mid W ],
g\in\G$.
Then $g_{0}$ satisfies the conditional moment restriction
%
\begin{equation}\label{moment}
m(W,g_{0}) = 0\mbox{,\qquad a.s.}
\end{equation}
Equivalently, we have $\E[ m^{2}(W,g_{0}) ] = 0$.

In this paper, for the purpose of robustness, any specific distribution
of $(Y,X,W)$ is not assumed, which we believe is more practical in
statistical and econometric applications. So a Bayesian analysis in the
standard sense is not applicable here since a proper likelihood for $g$
($g$ is a generic version of $g_{0}$) is not available.
Instead, we use a quasi-likelihood induced from the conditional moment
restriction (\ref{moment}).

Let $(Y_{1},X_{1},W_{1}),\ldots,(Y_{n},X_{n},W_{n})$ be i.i.d.
observations of $(Y,X,W)$.
Let $W^{n} = \{ W_{1},\ldots,W_{n} \}$ and $\D_{n} =\{
(Y_{1},X_{1},W_{1}),\ldots,(Y_{n},X_{n},W_{n}) \}$. By
(\ref{moment}),\break
a~plausible candidate of the quasi-likelihood would be
\[
p_{g}\bigl(W^{n}\bigr) = \exp\bigl\{ -(n/2) \En\bigl[
m^{2}(W_{i},g)\bigr] \bigr\},
\]
since $p_{g}(W^{n})$ is maximized at the true structural function $g_{0}$.
However, this $p_{g}(W^{n})$ is infeasible since $m(\cdot,g)$ is
unknown. Instead of using $p_{g}(W^{n})$, we replace $m(\cdot,g)$ by a
suitable estimate $\hat{m}(\cdot,g)$ and use the quasi-likelihood of
the form
\[
p_{g}(\D_{n}) = \exp\bigl\{ - (n/2) \En\bigl[
\hat{m}^{2}(W_{i},g)\bigr] \bigr\}.
\]
Below we use a wavelet series estimator of $m(\cdot,g)$.

The quasi-Bayesian analysis considered here uses this quasi-likelihood
as if it were a proper likelihood and puts priors on $g \in\G$. In this
paper, as in \cite{LJ11}, we shall use sieve priors (more precisely,
priors constructed on slowly growing sieves; \cite{LJ11}~indeed
considered another class of priors, see their supplementary material).
The basic idea is to construct a sequence of finite-dimensional sieves
(say, $\G_{n}$) that well approximates the parameter space $\G$ (i.e.,
each function in $\G$ is well approximated by some function in $\G_{n}$
as $n$ becomes large), and put priors concentrating on these sieves.
Each sieve space is a subset of a linear space spanned by some basis
functions. Hence the problem reduces to putting priors on the
coefficients on those basis functions. Such priors are typically called
``(finite dimensional) sieve priors'' (or ``series priors'') and have
been widely used in the nonparametric Bayesian and quasi-Bayesian
analysis (see, e.g., \cite{GGV00,S06,GV07}).

Let $\Pi_{n}$ be a so-constructed prior on $g \in\G$. Then,
formally, the posterior-like distribution of $g$ given $\D_{n}$ may be
defined by
%
\begin{equation}\label{qposterior}
\Pi_{n} ( dg \mid\D_{n}) = \frac{p_{g}(\D_{n}) \Pi_{n}(dg) }{\int
p_{g}(\D_{n}) \Pi_{n}(dg) },
\end{equation}
which we call ``quasi-posterior distribution.'' The quasi-posterior
distribution is not a proper posterior distribution in the strict
Bayesian sense since $p_{g}(\D_{n})$ is not a proper likelihood.
Nevertheless, $\Pi_{n} ( dg \mid\D_{n})$ is a proper distribution,
that is, $\int\Pi_{n}( dg \mid\D_{n}) = 1$. Similar to proper
posterior distributions,
contraction of the quasi-posterior distribution around $g_{0}$
intuitively means that it contains more and more accurate information
about the true structural function $g_{0}$ as the sample size increases.
Hence, as in proper posterior distributions, it is of fundamental
importance to study rates of contraction of quasi-posterior distributions.
Here we say that the quasi-posterior $\Pi_{n}(dg \mid\D_{n})$
contracts around $g_{0}$ at rate $\varepsilon_{n} \to0$ if $\Pi_{n}
( g\dvtx  \| g - g_{0} \| >\break \varepsilon_{n} \mid\D_{n} ) \stackrel
{P}{\to} 0$.

This quasi-posterior corresponds to what \cite{Z06b} called ``Gibbs
algorithm'' and what \cite{JT08} called ``Gibbs posterior.''
The framework of the quasi-posterior (Gibbs posterior) allows us a
flexibility since a stringent distributional assumption, such as
normality, on the data generating process is not required. Such a
framework widens a Bayesian approach to broad fields of statistical
problems.\footnote{Jiang and Tanner (\cite{JT08}, page 2211) remarked: ``This
framework of the Gibbs posterior has been overlooked by most
statisticians for a long time [$\cdots$] a foundation for
understanding the statistical behavior of the Gibbs posterior, which we
believe will open a productive new line of research.''}
Moreover, the following proposition gives an interesting interpretation
of the quasi-posterior.

\begin{proposition}
\label{prop0}
Let $\eta> 0$ be a fixed constant. Let $\Pi$ be a prior distribution
for $g$ defined on, say, the Borel $\sigma$-field of $C[0,1]$. Suppose
that the data $\D_{n}$ are fixed and the maps $g \mapsto\hat
{m}_{i}(W_{i},g)$ are measurable with respect to the Borel $\sigma
$-field of $C[0,1]$.
Then, the distribution
\[
\hat{\Pi}_{\eta} (dg) = \frac{\exp(- \eta\sum_{i=1}^{n} \hat
{m}^{2}(W_{i},g)) \Pi(dg)}{\int\exp(-\eta\sum_{i=1}^{n} \hat
{m}^{2}(W_{i},g)) \Pi(dg)}
\]
minimizes the empirical information complexity defined by
%
\begin{equation}\label{complexity}
\check{\Pi} \mapsto\int\sum_{i=1}^{n}
\hat{m}^{2}(W_{i},g) \check{\Pi}(dg) + \eta^{-1}
D_{\mathrm{KL}}( \check{\Pi} \,\|\, \Pi)
\end{equation}
over all distributions $\check{\Pi}$ absolutely continuous with
respect to $\Pi$. Here
\[
D_{\mathrm{KL}}( \check{\Pi} \,\|\, \Pi) = \int\check{\pi} \log\check{\pi}
\Pi(dg)\qquad \mbox{with } d\check{\Pi}/d \Pi= \check{\pi}
\]
is the Kullback--Leibler divergence from $\check{\Pi}$ to
$\Pi$.\vspace*{-2pt}
\end{proposition}
\begin{pf}
Immediate from \cite{Z06a}, Proposition 5.1.\vspace*{-2pt}
\end{pf}

The proposition shows that, given the data $\D_{n}$ and a prior $\Pi=
\Pi_{n}$ on $g$, the quasi-posterior $\Pi_{n}(dg \mid\D_{n})$ defined
in (\ref{qposterior}) is obtained as a minimizer of the empirical
information complexity defined by (\ref{complexity}) with $\eta= 1/2$.
This gives a rational to use $\Pi_{n}( dg \mid\D_{n})$ as a
quasi-posterior since, among all possible ``quasi-posteriors'', this
$\Pi_{n}( dg \mid\D_{n})$ optimally balances the average of the
natural loss function $g \mapsto\sum_{i=1}^{n} \hat{m}^{2}(W_{i},g)$
and its complexity (or deviation) relative to the initial prior
distribution measured by the Kullback--Leibler divergence. The scaling
constant (``temperature'') $\eta$ is typically treated as a fixed
constant (see, e.g., \mbox{\cite{Z06b,JT08}}). An alternative way is to
choose $\eta$ in a data-dependent manner, by, for example, cross
validation as mentioned in \cite{Z06b}. It is not difficult to see that
the theory below can be extended to the case where $\eta$ is even
random, as long as $\eta$ converges in probability to a fixed positive
constant. However, for the sake of simplicity, we take $\eta= 1/2$ as
a benchmark choice (note that as long as $\eta$ is a fixed positive
constant, the analysis can be reduced to the case with $\eta=1/2$ by
renormalization).

The quasi-posterior distribution provides point estimators of $g_{0}$.
A most natural estimator would be the estimator defined by
the posterior expectation (the expectation of the quasi-posterior
distribution), that is,
%
\begin{equation}
\label{qbe} \hat{g}_{\mathrm{QB}} = \cases{\displaystyle \int g \Pi_{n} ( dg \mid
\D_{n}), &\quad if the right integral exists,
\vspace*{2pt}\cr
0, &\quad otherwise,}
\end{equation}
where the integral $\int g \Pi_{n} ( dg \mid\D_{n})$ is understood
as pointwise.\eject

\begin{remark}
Quasi-Bayesian approaches (not necessarily in the present form) are
widely used and
there are several other attempts of making probabilistic interpretation
of such approaches.
See, for example, \cite{K02} where the ``limited information
likelihood'' is derived as the ``best'' (in a suitable sense)
approximation to the true likelihood function under a set of
moment restrictions and the Bayesian analysis with the limited
information likelihood is argued (\cite{LJ11} adapted this approach to
conditional moment restriction models), and
\cite{S05} where a version of the empirical likelihood is interpreted
in a Bayesian framework.
\end{remark}

\section{Main results}\label{sec3}
\label{secmain}

In this section, we study the asymptotic properties of the
quasi-posterior distribution and the quasi-Bayes estimator. In doing
so, we have to specify certain regularity properties, such as the
smoothness of $g_{0}$ and the degree of ill-posedness of the problem.
How to characterize the ``smoothness'' of $g_{0}$ is important since it
is related to how to put priors.
For this purpose, we find wavelet theory useful, and use sieve spaces
constructed by using wavelet bases.

\subsection{Posterior construction}\label{sec3.1}
\label{secpost}

To construct quasi-posterior distributions, we have to estimate
$m(\cdot,g)$ and construct a sequence of sieve spaces for $\G$ on which
priors concentrate. For the former purpose, we use a (wavelet) series
estimator of $m(\cdot,g)$, as in \cite{AC03} and \cite{CP11}. For the
latter purpose, we construct a sequence of sieve spaces formed by the
wavelet basis.

We begin with stating the parameter space for $g_{0}$ and the wavelet
basis used. We assume that the parameter space $\mathcal{G}$ is either
$(B^{s}_{\infty,\infty},\| \cdot\|_{s,\infty,\infty})$
(H\"{o}lder--Zygmund space) or $(B_{2,2}^{s},\| \cdot\|_{s,2,2})$
(Sobolev space), where $B^{s}_{p,q}$ is the Besov space of functions on
$[0,1]$ with parameter $(s,p,q)$ (the parameter $s$ generally
corresponds to ``smoothness;'' we add ``$s$'' on the parameter space,
$\mathcal {G}=\mathcal{G}^{s}$, to clarify its dependence on $s$). See
Appendix A.2 in the supplemental file \cite{K13} for the
definition of Besov spaces. We assume that $s > 1/2$, under which
$\mathcal{G}^{s} \subset C[0,1]$.

Fix (sufficiently large) $J_{0} \geq0$, and let $\{
\varphi_{J_{0}k}^{\inte} \}_{k=0}^{2^{J_{0}}-1} \cup\{
\psi_{jk}^{\inte}, j \geq J_{0}, k=0,\ldots,2^{j}-1 \}$ be
an $S$-regular Cohen--Daubechies--Vial (CDV) wavelet basis for
$L_{2}[0,1]$~\cite{CDV93}, where $S$ is a positive integer larger than
$s$.
See Appendix A.1 
in the supplemental file \cite{K13} for CDV wavelet bases.
For the notational convenience, we write $\phi_{1} = \varphi^{\inte
}_{J_{0},0},\phi_{2} = \varphi^{\inte}_{J_{0},1} ,\ldots, \phi
_{2^{J_{0}}} = \varphi^{\inte}_{J_{0},2^{J_{0}}-1}$, and $\phi
_{2^{j}+1} = \psi^{\inte}_{j,0}, \phi_{2^{j} + 2} = \psi^{\inte
}_{j,1}, \ldots, \phi_{2^{j+1}} = \psi^{\inte}_{j,2^{j}-1}$ for $j
\geq J_{0}$. Here and in what follows:
\[
\mbox{Take and fix an $S$-regular CDV wavelet basis of $\{ \phi_{l}, l \geq
1 \}$ with $S > s$,}
\]
and we keep this convention.
Let $V_{j}$ be the linear subspace of $L_{2}[0,1]$ spanned by $\{ \phi
_{1},\ldots,\phi_{2^{j}} \}$, and denote by $P_{j}$ the projection
operator onto $V_{j}$, that is, for any $g=\sum_{l=1}^{\infty} b_{l}
\phi_{l} \in L_{2}[0,1]$, $P_{j} g = \sum_{l=1}^{2^{j}} b_{l} \phi_{l}$.
In what follows, for any $J \in\mathbb{N}$, the notation $b^{J}$
means that it is a vector of dimension $2^{J}$. For example, $b^{J} =
(b_{1},\ldots,b_{2^{J}})^{T}$.

\begin{remark}[(Approximation property)]
\label{approx}
For either $g \in B^{s}_{\infty,\infty}$ or $B^{s}_{2,2}$, we have
$\| g - P_{J}g \|^{2} \leq C 2^{-2Js}$ for all $J \geq J_{0}$.
Here the constant $C$ depends only on $s$ and the corresponding Besov
norm of $g$.
\end{remark}

\begin{remark}
\label{bases}
The use of CDV wavelet bases is not crucial and one may use other
reasonable bases such as the Fourier and Hermite polynomial bases.
The theory below can be extended to such bases with some modifications.
However, CDV wavelet bases are particularly well suited to approximate
(not necessarily periodic) smooth functions, which is the reason why we
use here CDV wavelet bases. On the other hand, for example, the Fourier
basis is only appropriate to approximate periodic functions and it is
often not natural to assume that the structural function $g_{0}$ is periodic.
\end{remark}

We shall now move to the posterior construction.
For $J \geq J_{0}$, define the $2^{J}$-dimensional vector of functions
$\phi^{J}(w)$ by
\[
\phi^{J}(w) = \bigl( \phi_{1}(w),\ldots,
\phi_{2^{J}}(w)\bigr)^{T}.
\]
Let $J_{n} \geq J_{0}$ be a sequence of positive integers such that
$J_{n} \to\infty$ and $2^{J_{n}} = o(n)$. Then a wavelet series
estimator of $m(\cdot,g)$ is defined as
\[
\hat{m}(w,g) = \phi^{J_{n}} (w)^{T} \bigl(\En\bigl[ \phi
^{J_{n}}(W_{i})^{\otimes2}\bigr]\bigr)^{-1} \En
\bigl[ \phi^{J_{n}} (W_{i}) \bigl(Y_{i} -
g(X_{i})\bigr)\bigr],
\]
where we replace the inverse matrix by the generalized inverse if the
former does not exist; the probability of such an event converges to
zero as $n \to\infty$ under the assumptions below. We use this
wavelet series estimator throughout the analysis.

For the same $J_{n}$, we shall take $V_{J_{n}} = \spn\{
\phi_{1},\ldots,\phi_{2^{J_{n}}} \}$ as a sieve space for~$\G^{s}$. We
consider priors $\Pi_{n}$ that concentrate on $V_{J_{n}}$, that is,
$\Pi_{n}(V_{J_{n}}) = 1$. Formally, we think of that priors on $g$ are
defined on the Borel $\sigma$-field of $C[0,1]$ (hence the
quasi-posterior $\Pi_{n}(dg \mid\D_{n})$ is understood to be defined on
the Borel $\sigma$-field of $C[0,1]$, which is possible since the map
$g \mapsto p_{g}(\D_{n})$ is continuous on $C[0,1]$). Since the map
$b^{J_{n}} = (b_{1},\ldots,b_{2^{J_{n}}})^{T}
\mapsto\sum_{l=1}^{2^{J_{n}}} b_{l} \phi_{l}, \R^{2^{J_{n}}} \to
C[0,1]$, is homeomorphic from $\R^{2^{J_{n}}}$ onto $V_{J_{n}}$,
putting priors on $g \in V_{J_{n}}$ is equivalent to putting priors on
$b^{J_{n}} \in\R^{2^{J_{n}}}$ (the latter are of course defined on the
Borel $\sigma$-field of $\R^{2^{J_{n}}}$). Practically, priors on $g
\in V_{J_{n}}$ are induced from priors on $b^{J_{n}}
\in\R^{2^{J_{n}}}$. For the later purpose, it is useful to determine
the correspondence between priors for these two parameterizations.
Unless otherwise stated, we follow the convention of the notation such
that
\[
\mbox{$\tilde{\Pi}_{n}$: a prior on $b^{J_{n}} \in\R^{2^{J_{n}}}$
$\leftrightarrow$ $\Pi_{n}$: the induced prior on $g \in V_{J_{n}}$.}
\]
We shall call $\tilde{\Pi}_{n}$ a generating prior, and $\Pi_{n}$
the induced prior.\eject

Correspondingly, the quasi-posterior for $b^{J_{n}}$ is defined. With a
slight abuse of notation, for $g = \sum_{l=1}^{2^{J_{n}}} b_{l} \phi
_{l}$, we write\vspace*{1pt} $\hat{m}(w,b^{J_{n}}) = \hat{m}(w,g)$, and take
$p_{b^{J_{n}}} (\D_{n}) = \exp\{ -(n/2) \En[\hat
{m}^{2}(W_{i},b^{J_{n}})] \}$ as a quasi-likelihood for $b^{J_{n}}$.
Note that\vspace*{1pt} in this particular setting, the log quasi-likelihood is
quadratic in $b^{J_{n}}$.
Let $\tilde{\Pi}_{n}( db^{J_{n}} \mid\D_{n})$ denote the resulting
quasi-posterior distribution for $b^{J_{n}}$:
%
\begin{equation}\label{qposterior-b}
\tilde{\Pi}_{n}\bigl( db^{J_{n}} \mid\D_{n}\bigr) =
\frac{p_{b^{J_{n}}} (\D
_{n}) \tilde{\Pi}_{n}(db^{J_{n}})}{\int p_{b^{J_{n}}} (\D_{n})
\tilde{\Pi}_{n}(db^{J_{n}})}.
\end{equation}

For the quasi-Bayes estimator $\hat{g}_{\mathrm{QB}}$ defined by (\ref{qbe}),
since for every $x \in[0,1]$, the map $g \mapsto g(x)$ is continuous
on $C[0,1]$, and conditional on $\D_{n}$ the quasi-posterior $\Pi
_{n}(dg \mid\D_{n})$ is a Borel probability measure on $C[0,1]$, the
integral $\int g(x) \Pi_{n}(dg \mid\D_{n})$ exists as soon as $\int
| g(x) | \Pi_{n}(dg \mid\D_{n}) < \infty$.
Furthermore, $\hat{g}_{\mathrm{QB}}$ can be computed by using the relation
\[
\int g(x) \Pi_{n}(dg \mid\D_{n}) = \phi^{J_{n}}(x)^{T}
\biggl[ \int b^{J_{n}} \tilde{\Pi}_{n}\bigl(db^{J_{n}}
\mid\D_{n}\bigr) \biggr]
\]
as soon as the integral on the right-hand side exists.
Hence, practically, it is sufficient to compute the expectation of
$\tilde{\Pi}_{n}(db^{J_{n}} \mid\D_{n})$.

\begin{remark}
The use of the same wavelet basis to estimate $m(\cdot,g)$ and to
construct a sequence of sieve spaces for $\G^{s}$ is not essential and
can be relaxed. Suppose that we have another CDV wavelet basis $\{
\tilde{\phi}_{l} \}$ for $L_{2}[0,1]$ and use this basis to estimate
$m(\cdot,g)$. Then, all the results below apply by simply replacing
$\phi_{l}(W_{i})$ by $\tilde{\phi}_{l}(W_{i})$.
To keep the notation simple, we use the same wavelet basis.

However, the use of the same resolution level $J_{n}$ is essential (at
least at the proof level) in establishing the asymptotic properties of
the quasi-posterior distribution. It may be a technical artifact, but
we do not extend the theory in this direction since there is no clear
theoretical benefit to do so (note that in the purely frequentist
estimation case, \cite{CP11} allowed for using different cut-off levels
for approximating $m(\cdot,g)$ and $g(\cdot)$).
\end{remark}

\subsection{Basic assumptions}\label{sec3.2}
\label{secassumption}

We state some basic assumptions. We do not state here assumptions on
priors, which will be stated in the theorems below.
In what follows, let $C_{1} > 1$ be a sufficiently large constant.
%
\begin{assumption}
\label{A1}
(i) $(X,W)$ has a joint density $f_{X,W}(x,w)$ on $[0,1]^{2}$
satisfying that $f_{X,W}(x,w) \leq C_{1}, \forall x,w \in[0,1]$.
(ii) $\sup_{w \in[0,1]}\E[ U^{2} \mid W=w ] \leq C_{1}$ where $U=Y-g_{0}(X)$.
(iii) $s_{\min}(\E[ \phi^{J}(W)^{\otimes2} ]) \geq C_{1}^{-1},
\forall J \geq J_{0}$.
\end{assumption}
Assumption \ref{A1} is a usual restriction in the literature, up to
minor differences (see \cite{HH05,H11a}). Denote by $f_{X} (x)$ and
$f_{W}(w)$ the marginal densities of $X$ and $W$, respectively, that
is, $f_{X}(x) = \int f_{X,W}(x,w) \,dw$ and $f_{W}(w) = \int
f_{X,W}(x,w) \,dx$. Then Assumption \ref{A1}(i) implies that $f_{X}(x)
\leq C_{1}, \forall x \in[0,1]$ and $f_{W}(w) \leq C_{1}, \forall w
\in[0,1]$. A primitive regularity condition that guarantees Assumption
\ref{A1}(iii) is that $f_{W}(w) \geq C^{-1}_{1}$ for all $w \in[0,1]$.
To see this, for $\alpha^{J} \in\R^{2^{J}}$ with $\| \alpha^{J}
\|_{\ell ^{2}}=1$, we have
%
\begin{eqnarray}\label{lowereigen}
\bigl(\alpha^{J}\bigr)^{T}\E\bigl[ \phi^{J}(W)^{\otimes2}
\bigr]\alpha^{J} &=& \int_{0}^{1} \bigl(
\phi^{J}(w)^{T}\alpha^{J}\bigr)^{2}
f_{W}(w) \,dw \nonumber\\
&\geq& C_{1}^{-1} \int
_{0}^{1}\bigl(\phi^{J}(w)^{T}
\alpha^{J}\bigr)^{2} \,dw
\nonumber\\[-8pt]\\[-8pt]
& = &C_{1}^{-1} \bigl(\alpha^{J}
\bigr)^{T} \biggl[ \int_{0}^{1}
\phi^{J}(w)\phi^{J}(w)^{T} \,dw \biggr]
\alpha^{J} \nonumber\\
&=& C_{1}^{-1}\bigl\| \alpha^{J}
\bigr\|_{\ell
^{2}}^{2}= C_{1}^{-1},\nonumber
\end{eqnarray}
where we have used the fact that $\{ \phi_{l} \}$ is orthonormal in
$L_{2}[0,1]$.

For identification of $g_{0}$, we assume:
%
\begin{assumption}
\label{A2}
The linear operator $K\dvtx  L_{2}[0,1] \to L_{2}[0,1]$ is injective.
\end{assumption}

For smoothness of $g_{0}$, as mentioned before, we assume:
%
\begin{assumption}
\label{A3}
$\exists s > 1/2$, $g_{0} \in\G^{s}$, where $\G^{s}$ is either
$B^{s}_{\infty,\infty}$ or $B^{s}_{2,2}$.
\end{assumption}

The identification condition (Assumption \ref{A2}) is equivalent to
the ``completeness'' of the conditional distribution of $X$ conditional
on $W$ \cite{NP03}.
We refer the reader to \cite{ST06,D11} and \cite{HS11} for discussion
on the completeness condition. We should note that restricting the
domain of $K$ to a ``small'' set, such as a Sobolev ball, would
substantially relax Assumption \ref{A2}, which however requires a
different analysis. For the sake of simplicity,
we assume the injectivity of $K$ on the full domain.

As discussed in the \hyperref[sec1]{Introduction}, solving (\ref{npiv2}) is
an ill-posed inverse problem. Thus, the statistical difficulty of
estimating $g_{0}$ depends on the difficulty of continuously inverting
$K$, which is usually referred to as ``ill-posedness'' of the inverse
problem~(\ref{npiv2}). Typically, the ill-posedness is characterized by
the decay rate of $\kappa_{l} \to0$ ($\kappa_{l}$ is the $l$th largest
singular value of $K$), which is plausible if $K$ were known and the
singular value decomposition of $K$ were used (see \cite{C08}).
However, here, $K$ is unknown and the known wavelet basis $\{ \phi_{l}
\}$ is used instead of the singular value system. Thus, it is suitable
to quantify the ill-posedness using the wavelet basis~$\{ \phi_{l} \}$.
To this end, define
\[
\tau_{J} = s_{\min}\bigl(\E\bigl[ \phi^{J}(W)
\phi^{J}(X)^{T}\bigr]\bigr) = s_{\min
} \bigl( \bigl(
\langle\phi_{l},K \phi_{m} \rangle\bigr)_{1 \leq l,m \leq
2^{J}}
\bigr),\qquad J \geq J_{0}.
\]

This quantity corresponds to (the reciprocal of) what is called ``sieve
measure of ill-posedness'' in the literature \cite{BCK07,H11a}. We
at\vadjust{\goodbreak}
least have to assume that $\tau_{J} > 0$ for all $J \geq J_{0}$. Note
however that
\begin{eqnarray*}
\tau_{J} &=& s_{\min} \bigl( \bigl(\langle
\phi_{l},K \phi_{m} \rangle\bigr)_{1 \leq l,m \leq2^{J}} \bigr)
\\
&=&\min_{g \in V_{J}, \| g \|=1} \bigl\| \bigl(\langle\phi_{l},K g \rangle
\bigr)_{1
\leq l \leq2^{J}} \bigr\|_{\ell^{2}}
\\
&\leq&\min_{g \in V_{J}, \| g \|=1} \| K g \| \qquad\mbox{(Bessel's
inequality)}
\\
&\leq&\kappa_{2^{J}}\qquad\mbox{(Courant--Fischer--Weyl's
minimax principle)}
\end{eqnarray*}
by which, necessarily, $\tau_{J} \to0$ as $J \to\infty$.
For this quantity, we assume:
%
\begin{assumption}
\label{A4}
(i) (Mildly ill-posed case) $\exists r > 0$, $\tau_{J} \geq C^{-1}_{1}
2^{-Jr}, \forall J \geq J_{0}$ or (severely ill-posed case) $\exists
c > 0$, $\tau_{J} \geq C^{-1}_{1} \exp( - c 2^{J}), \forall J \geq J_{0}$;

(ii)
\begin{eqnarray*}
&&\bigl\| \E\bigl[ \phi^{J}(W) (g_{0}-P_{J}g_{0})
(X) \bigr] \bigr\|_{\ell^{2}} \bigl(= \bigl\| \bigl( \bigl\langle\phi_{l},
K(g_{0} - P_{J}g_{0}) \bigr\rangle
\bigr)_{l=1}^{2^{J}} \bigr\| _{\ell^{2}}\bigr)
\\
&&\qquad
\leq C_{1} \tau_{J} \| g_{0}-P_{J}g_{0}
\|\qquad \forall J \geq J_{0}.
\end{eqnarray*}
\end{assumption}

Assumption \ref{A4}(i) lower bounds $\tau_{J}$ as $J \to\infty$,
thereby quantifies the ill-posedness. We cover both the ``mildly
ill-posed'' and ``severely ill-posed'' cases (this definition of mild
ill-posedness and severe ill-posedness is due to \mbox{\cite{H11b,H11a}}).
The severely ill-posed case happens, for example, when the joint
density $f_{X,W}(x,w)$ is analytic (see \cite{K99}, Theorem 15.20).

Assumption \ref{A4}(ii) is a ``stability'' condition about the bias
$g_{0} - P_{J} g_{0}$, which states that $K(g_{0} - P_{J}g_{0})$ is
sufficiently ``small'' relative to $g_{0} - P_{J}g_{0}$. Note that in
the (ideal) case in which, for example, $K$ is self-adjoint and $\{
\phi_{l} \}$ is the eigen-basis of $K$, $\langle\phi_{l}, K(g_{0} -
P_{J}g_{0}) \rangle= 0$ for all $l=1,\ldots,2^{J}$, in which case
Assumption~\ref{A4}(ii) is trivially satisfied. Assumption \ref{A4}(ii)
allows more general situations in which $K$ may not be self-adjoint and
$\{ \phi_{l} \}$ may not be the eigen-basis of $K$ by allowing for a
certain ``slack.'' This assumption, although looks technical, is common
in the study of rates of convergence in estimation of the structural
function~$g_{0}$. Indeed, essentially similar conditions have appeared
in the past literature such as \mbox{\cite{BCK07,CR11,H11a}}. For example,
\cite{BCK07}, Assumption 6, essentially states (in our notation) that
$\| K(g_{0} - P_{J}g_{0}) \| \leq C_{1} \tau_{J} \| g_{0} - P_{J}g_{0}
\|$, which implies our Assumption \ref{A4}(ii) since $\| (
\langle\phi_{l}, K(g_{0} - P_{J}g_{0}) \rangle)_{l=1}^{2^{J}}
\|_{\ell^{2}} \leq\| K(g_{0} - P_{J}g_{0}) \|$ (Bessel's inequality).

\begin{remark}
\label{minimax}
For given values of $C_{1} > 1,M > 0,r>0,c>0$ and $s>1/2$,
let $\mathcal{F} = \mathcal{F}(C_{1},M,r,c,s)$ denote the set of all
distributions of $(Y,X,W)$ satisfying Assumptions \ref{A1}--\ref{A4}
with $\| g_{0} \|_{s,\infty,\infty} \leq M$ in case of $\G^{s} =
B^{s}_{\infty,\infty}$ and $\| g_{0} \|_{s,2,2} \leq M$ in case of
$\G^{s} = B^{s}_{2,2}$.
By \cite{HH05,CR11}, it is shown that
the minimax rate of convergence (in $\| \cdot\|$) of estimation of
$g_{0}$ over this distribution class $\mathcal{F}$ is
$n^{-s/(2r+2s+1)}$ in the mildly ill-posed case (where $\tau_{J} \geq
C_{1}^{-1} 2^{-Jr}$) and
$(\log n)^{-s}$ in the severely ill-posed case [where $\tau_{J} \geq
C_{1}^{-1} \exp(-c2^{J})$] as the sample size $n \to\infty$ (the
assumption on the conditional second moment of $U$ given $W$ is not
binding;\vspace*{1pt} that is, replacing Assumption \ref{A1}(ii) by a stronger one, such as
$\sup_{w \in[0,1]} \E[ | U |^{2+\epsilon} \mid W=w] \leq C_{1}$ for
some $\epsilon> 0$ determined outside the class of distributions, does
not alter these minimax rates).

By Theorem 2.5 of \cite{GGV00}, it is readily seen that these rates
are the fastest possible rates of contraction of (general)
quasi-posterior distributions in this setting. More formally, we can
state the following assertion:\vspace*{9pt}

\emph{Let $\Pi_{n}( dg \mid\D_{n})$ be the quasi-posterior distribution
defined on, say, the Borel $\sigma$-field of $C[0,1]$, constructed
from putting a suitable prior on $g$ to the quasi-likelihood $p_{g}(\D
_{n})$} (\textit{the prior here needs not be a sieve prior}).
\textit{Suppose now that for some $\varepsilon_{n} \to0, \sup_{F \in
\mathcal{F}} \E_{F} [ \Pi_{n}( g\dvtx  \| g - g_{0} \| > \varepsilon
_{n} \mid\D_{n}) ] \to0$.
Then there exists a point estimator that converges} (\textit{in probability})
\textit{at least as fast as $\varepsilon_{n}$ uniformly in $F
\in\mathcal{F}$.}\vspace*{9pt}

The proof is just a small modification of that of Theorem 2.5 in
\cite{GGV00} and hence omitted. Importantly, the quasi-posterior cannot
contract at a rate faster than the optimal rate of convergence for
point estimators (\cite{GGV00}, page 507, lines 19--20). Hence, in
the minimax sense, the fastest possible rate of contraction of the
quasi-posterior distribution $\Pi_{n}(dg \mid\D_{n})$ is
$n^{-s/(2r+2s+1)}$ in the mildly ill-posed case and $(\log n)^{-s}$ in
the severely ill-posed case (Proposition \ref{prop1} in Section \ref
{secexamples} ahead shows that these rates are indeed attainable
for suitable sieve priors).
\end{remark}

\subsection{Main results: General theorems}\label{sec3.3}

This section presents general theorems on contraction rates and
asymptotic normality for the quasi-posterior distribution as well as
convergence rates for the quasi-Bayes estimator.
In what follows, let $(Y_{1},X_{1},W_{1}),\ldots,(Y_{n},X_{n},W_{n})$
be i.i.d. observations of $(Y,X,W)$.
Denote by $b_{0}^{J} = (b_{01},\ldots,b_{0,2^{J}})^{T}$ the vector of
the first $2^{J}$ generalized Fourier coefficients of~$g_{0}$, that is,
$b_{0l} = \int\phi_{l} g_{0}$. Let $\| \cdot\|_{\TV}$ denote the
total variation norm between two distributions.

\begin{theorem}
\label{thm1}
Suppose that Assumptions \ref{A1}--\ref{A4} are satisfied. Take
$J_{n}$ in such a way that $J_{n} \to\infty$ and $J_{n}2^{J_{n}}/n =
o(\tau_{J_{n}}^{2})$. Let $\epsilon_{n}$ be a sequence of positive
constants such that $\epsilon_{n} \to0$ and $n \epsilon_{n}^{2}
\gtrsim2^{J_{n}}$.
Suppose that generating priors $\tilde{\Pi}_{n}$ has densities
$\tilde{\pi}_{n}$ on $\R^{2^{J_{n}}}$ and satisfy the following conditions:
\begin{longlist}[(P2)]
\item[(P1)] (Small ball condition). There exists a constant $C>0$ such
that for all $n$ sufficiently large, $\tilde{\Pi}_{n} (b^{J_{n}}\dvtx  \|
b^{J_{n}} - b_{0}^{J_{n}} \|_{\ell^{2}} \leq\epsilon_{n}) \geq
e^{-Cn\epsilon_{n}^{2}}$.
\item[(P2)] (Prior flatness condition). Let $\gamma_{n} = 2^{-J_{n}s} +
\tau_{J_{n}}^{-1} \epsilon_{n}$. There exists a sequence of constants
$L_{n} \to\infty$ sufficiently slowly such that for all $n$
sufficiently large, $\tilde{\pi}_{n}(b^{J_{n}})$ is positive for all
$\| b^{J_{n}}-b_{0}^{J_{n}} \|_{\ell^{2}} \leq L_{n} \gamma_{n}$, and
\[
\sup_{\| b^{J_{n}} \|_{\ell^{2}} \leq L_{n} \gamma_{n}, \| \tilde
{b}^{J_{n}} \|_{\ell^{2}} \leq L_{n} \gamma_{n}} \biggl\llvert\frac
{\tilde{\pi}_{n}(b_{0}^{J_{n}}+b^{J_{n}})}{\tilde{\pi
}_{n}(b_{0}^{J_{n}}+\tilde{b}^{J_{n}})} - 1 \biggr\rrvert
\to0.
\]
\end{longlist}
Then for every sequence $M_{n} \to\infty$, we have
%
\begin{equation}\label{rate}
\tilde{\Pi}_{n} \bigl\{ b^{J_{n}}\dvtx  \bigl\| b^{J_{n}} -
b^{J_{n}}_{0} \bigr\| _{\ell^{2}} > M_{n}
\bigl(2^{-J_{n}s} + \tau_{J_{n}}^{-1} \sqrt
{2^{J_{n}}/n}\bigr) \mid\D_{n} \bigr\} \stackrel{P} {\to} 0.
\end{equation}
Furthermore, assume that $J_{n}2^{3J_{n}}/n = o(\tau_{J_{n}}^{2})$.
Then we have
\[
\bigl\| \tilde{\Pi}_{n} (\cdot\mid\D_{n})- N\bigl(
\hat{b}^{J_{n}}, n^{-1} \Phi^{-1}_{WX}
\Phi_{WW} \Phi_{XW}^{-1}\bigr) (\cdot)
\bigr\|_{\TV} \stackrel{P} {\to} 0,
\]
where $\Phi_{WX}:= \E[ \phi^{J_{n}}(W)\phi^{J_{n}}(X)^{T} ], \Phi
_{XW}:=\Phi_{WX}^{T}, \Phi_{WW}:=\E[ \phi^{J_{n}}(W)^{\otimes2}]$,
and where $\hat{b}^{J_{n}}$ is a ``maximum quasi-likelihood
estimator'' of $b_{0}^{J_{n}}$, that is,
%
\begin{equation}\label{BvM}
\hat{b}^{J_{n}} \in\arg\max_{b^{J_{n}} \in\R^{2^{J_{n}}}} p_{b^{J_{n}}}(
\D_{n}).
\end{equation}
\end{theorem}

\begin{pf}
See Section \ref{secprf1}.
\end{pf}

\begin{remark}
The condition $J_{n}2^{J_{n}}/n = o(\tau_{J_{n}}^{2})$ appears
essentially because the operator $K$ is unknown. In our setup, this
results in estimating the matrix $\E[\phi^{J_{n}}(W) \phi
^{J_{n}}(X)^{T} ]$ by its empirical counterpart $\En[\phi
^{J_{n}}(W_{i}) \phi^{J_{n}}(X_{i})^{T} ]$. In the proof, we have to
suitably lower bound the minimum singular value of $\En[\phi
^{J_{n}}(W_{i}) \phi^{J_{n}}(X_{i})^{T}]$, denoted by $\hat{\tau
}_{J_{n}}$, which is an empirical counterpart of the sieve measure of
ill-posedness $\tau_{J_{n}}$. By Lemma \ref{lemB1}, we have $\hat
{\tau}_{J_{n}} = \tau_{J_{n}}-O_{P}(\sqrt{J_{n} 2^{J_{n}}/n})$, so
that to make the estimation effect in $\hat{\tau}_{J_{n}}$
negligible, we need $J_{n} 2^{J_{n}}/n = o(\tau_{J_{n}}^{2})$.
\end{remark}

\begin{remark}
Theorem \ref{thm1} is abstract in the sense that it only gives
conditions (P1) and (P2) on priors for which (\ref{rate}) and (\ref
{BvM}) hold. For specific priors, we have to check these conditions
with possible $J_{n}$, which will be done in Section \ref{secexamples}.
\end{remark}

Since for $g = \sum_{l=1}^{2^{J_{n}}} b_{l} \phi_{l}$, $\| g - g_{0}
\|^{2} = \| g - P_{J_{n}}g_{0} \|^{2} + \| g_{0} - P_{J_{n}} g_{0} \|
^{2} \lesssim\| b^{J_{n}} - b_{0}^{J_{n}} \|_{\ell^{2}}^{2} +
2^{-2J_{n}s}$, part (\ref{rate}) of Theorem \ref{thm1} leads to that
for every sequence $M_{n} \to\infty$, we have
\[
\Pi_{n} \bigl\{ g\dvtx  \| g -g_{0} \| > M_{n}
\bigl(2^{-J_{n}s} + \tau_{J_{n}}^{-1} \sqrt{2^{J_{n}}/n}
\bigr) \mid\D_{n} \bigr\} \stackrel{P} {\to} 0,
\]
which means that the rate of contraction of the quasi-posterior
distribution $\Pi_{n}( dg \mid\D_{n})$ is $\max\{ 2^{-J_{n}s}, \tau
_{J_{n}}^{-1}\sqrt{2^{J_{n}}/n} \}$.\footnote{We have ignored the
appearance of $M_{n} \to\infty$, which can be arbitrarily slow. A
version in which $M_{n}$ is replaced by a large fixed constant $M > 0$
is presented in Theorem \ref{thm3}.}
In many examples, for given $J_{n} \to\infty$ with $J_{n}2^{J_{n}}/n
= o(\tau_{J_{n}}^{2})$, condition (P1) is satisfied with $\epsilon_{n}
\sim\sqrt{2^{J_{n}}(\log n) /n}$.
Taking $J_{n}$ in such a way that [with some constant $c' < 1/(2c)$ in
the severely ill-posed case]
%
\begin{equation}
\label{optimal-J} \cases{ 2^{J_{n}} \sim n^{1/(2r+2s+1)}, &\quad in the mildly
ill-posed case,
\vspace*{2pt}\cr
\displaystyle \lim_{n \to\infty} \bigl( 2^{J_{n}}/
\bigl(c' \log n\bigr) \bigr)=1, &\quad in the severely ill-posed case,}
\end{equation}
under which the optimal contraction rate is attained, $\gamma_{n}$ in
condition (P2) is
%
\begin{equation}
\label{radius} \gamma_{n} \sim\cases{ n^{-s/(2r+2s+1)} (\log
n)^{1/2}, &\quad in the mildly ill-posed case,
\cr
(\log n)^{-s}, &\quad
in the severely ill-posed case.}
\end{equation}
So condition (P2) states that, to attain the optimal contraction rate
(and the Bernstein--von Mises type result), the prior density $\tilde
{\pi}_{n}$ should be sufficiently ``flat'' in a ball with center
$b_{0}^{J_{n}}$ and radius of order (\ref{radius}). Some specific
priors leading to the optimal contraction rate will be given in
Section \ref{secexamples}.

As noted before, in many examples, for given $J_{n} \to\infty$ with
$J_{n}2^{J_{n}}/n = o(\tau_{J_{n}}^{2})$, condition\vspace*{1pt} (P1) is satisfied
with $\epsilon_{n} \sim\sqrt{2^{J_{n}}(\log n)/n}$.
Inspection of the proof shows that, without condition (P2), this already
leads to contraction rate $\max\{ 2^{-J_{n}s}, \tau_{J_{n}}^{-1}
\sqrt{2^{J_{n}}(\log n)/n} \}$, which, in the mildly ill-posed case,
reduces to $(n/\log n)^{-s/(2r+2s+1)}$ by taking $2^{J_{n}} \sim
(n/\log n)^{1/(2r+2s+1)}$. However, this rate is not fully satisfactory
because of the appearance of the log term. Condition (P2) is used to get
rid of the log term.

The small ball condition (P1) is standard in nonparametric Bayesian
statistics and analogous to condition (2.4) in \cite{GGV00}. It is,
however, stated in
\cite{GGV00}, pages 505--506, that their Theorem 2.1 is not sharp
enough when
priors constructed on a sequence of finite-dimensional sieves are used,
and the more sophisticated condition (2.9) is devised in their Theorem 2.4
(see also the proof of their Theorem 4.5). However, a version of their
condition (2.9) is not clear to work in our problem, because the effect
of the random matrix $\En[ \phi^{J_{n}}(W_{i}) \phi
^{J_{n}}(X_{i})^{T} ]$ has to be suitably controlled. Instead, we
devise condition (P2) to obtain sharper contraction rates.

Under a further integrability condition about $U = Y-g_{0}(X)$, $M_{n}
\to\infty$ in (\ref{rate}) can be replaced by a large fixed constant $M$.
%
\begin{theorem}
\label{thm3}
Suppose that all the conditions that guarantee (\ref{rate}) in Theorem
\ref{thm1} are satisfied. Furthermore, assume that $\sup_{w \in
[0,1]}\E[ U^{2} 1( | U | > \lambda) \mid W=w] \to0$ as $\lambda\to
\infty$ where $U=Y-g_{0}(X)$. Then there exists a constant $M > 0$
such that
%
\begin{equation}\label{rate2}
\tilde{\Pi}_{n} \bigl\{ b^{J_{n}}\dvtx  \bigl\| b^{J_{n}} -
b^{J_{n}}_{0} \bigr\| _{\ell^{2}} > M \bigl(2^{-J_{n}s} +
\tau_{J_{n}}^{-1}\sqrt{2^{J_{n}}/n}\bigr) \mid
\D_{n} \bigr\} \stackrel{P} {\to} 0.
\end{equation}
\end{theorem}

\begin{pf}
See Section \ref{secprf3}.
\end{pf}\eject

The proof consists in establishing a concentration property of the
random variable $\| \En[ \phi^{J_{n}}(W_{i}) U_{i} ] \|_{\ell^{2}}$,
which uses a truncation argument and Talagrand's \cite{T96}
concentration inequality. A sufficient condition that guarantees
that
\[
\sup_{w \in[0,1]}\E\bigl[ U^{2} 1\bigl(|U| > \lambda\bigr) \mid W=w\bigr] \to0
\]
as\vspace*{1pt} $\lambda\to\infty$ is that $\exists\epsilon> 0$, $\sup_{w \in [0,1]}
\E[ | U |^{2+\epsilon} \mid W=w] < \infty$. The additional condition in
Theorem \ref{thm3} is a uniform integrability condition and stronger
than Assumption \ref{A1}(ii). To see this, note that $U$ is distributed
as $F^{-1}_{U \mid W}(\mathcal{U} \mid W)$ where\vspace*{1pt} $F^{-1}_{U \mid W}(u
\mid w)$ is the conditional quantile function of $U$ given $W=w$, and
$\mathcal{U}$ is a uniform random variable on $(0,1)$ independent of
$W$. Think of $U_{w}(u) = F^{-1}_{U \mid W}(u \mid w), w \in[0,1]$ as a
stochastic process defined on the probability space $((0,1), \mu)$ with
$\mu$ Lebesgue measure on $(0,1)$. Then the\vspace*{1pt} condition $\sup_{w
\in[0,1]}\E[ U^{2} 1( | U | > \lambda) \mid W=w] \to0$ as
$\lambda\to\infty$ states exactly the uniform integrability of
$(U_{w})_{w \in[0,1]}$.

The second part of Theorem \ref{thm1} states a Bernstein--von Mises type
result for the quasi-posterior distribution $\tilde{\Pi}_{n}(db^{J_{n}}
\mid\D_{n})$, which states that the quasi-posterior distribution is
approximated by the normal distribution centered at $\hat{b}^{J_{n}}$,
which is often referred to as the ``sieve minimum distance estimator''
and is a benchmark frequentist estimator for these types of models.
Note that, neglecting the bias, $\hat{b}^{J_{n}}$~is approximated as
$b^{J_{n}}_{0} + \Phi_{WX}^{-1} \En[ \phi^{J_{n}}(W_{i}) U_{i} ]$, but
the covariance matrix of the term $\Phi_{WX}^{-1} \sqrt{n}
\En[\phi^{J_{n}}(W_{i})U_{i}]$ is generally different from
$\Phi_{WX}^{-1}\Phi_{WW}\Phi_{XW}^{-1}$ (which is the reason why we
added ``type''). This is a generic nature of quasi-posterior
distributions. Even for finite-dimensional models, generally, the
covariance matrix of the centering variable does not coincide with
that of the normal distribution approximating the quasi-posterior
distribution (see \cite{CH03}).

Finally, we consider the convergence rate of the quasi-Bayes estimator
$\hat{g}_{\mathrm{QB}}$ of $g_{0}$ defined by (\ref{qbe}).
%
\begin{theorem}
\label{thm4}
Suppose that all the conditions of Theorem \ref{thm3} are satisfied.
Let $\hat{g}_{\mathrm{QB}}$ be the quasi-Bayes estimator defined by (\ref
{qbe}). Then $\Prob\{ \D_{n}\dvtx\break  \int| g(x) | \Pi_{n}(dg \mid\D
_{n}) < \infty, \forall x \in[0,1] \} \to1$, and there exists a
constant $D > 0$ such that for every sequence $M_{n} \to\infty$,
%
\begin{equation}\label{convrate}
\Prob\bigl[ \| \hat{g}_{\mathrm{QB}} - g_{0} \|
\leq D \max\bigl\{ 2^{-J_{n}s}, \tau_{J_{n}}^{-1}
\sqrt{2^{J_{n}}/n}, \tau_{J_{n}}^{-1}\epsilon_{n}
\varrho_{n}M_{n} \bigr\} \bigr] \to1,
\end{equation}
where
\[
\varrho_{n}:= \sup_{\| b^{J_{n}} \|_{\ell^{2}} \leq L_{n} \gamma
_{n}, \| \tilde{b}^{J_{n}} \|_{\ell^{2}} \leq L_{n} \gamma_{n}} \biggl
\llvert
\frac{\tilde{\pi}_{n}(b_{0}^{J_{n}}+b^{J_{n}})}{\tilde{\pi
}_{n}(b_{0}^{J_{n}}+\tilde{b}^{J_{n}})} - 1 \biggr\rrvert,
\]
and where $\epsilon_{n}, \gamma_{n}$ and $L_{n}$ are given in the
statement of Theorem \ref{thm1}.
\end{theorem}

\begin{pf}
See Appendix C 
in the supplemental file \cite{K13}.
\end{pf}

Theorem \ref{thm4} is \emph{not} directly deduced from Theorem \ref
{thm1}. Indeed, $\| g - g_{0} \|$ may be unbounded on the support of
$\Pi_{n}$ since the support of $\Pi_{n}$ may be unbounded in $\|
\cdot\|$, and hence
the argument in \cite{GGV00}, pages 506--507, cannot apply (in \cite
{GGV00}, a typical distance to measure the goodness of a point
estimator is the Hellinger distance and uniformly bounded).
Hence, additional work is needed to prove Theorem \ref{thm4}.

The convergence rate of the quasi-Bayes estimator is determined by the
three terms: $2^{-J_{n}s}, \tau_{J_{n}}^{-1} \sqrt{2^{J_{n}}/n}$, and
$ \tau_{J_{n}}^{-1} \epsilon_{n} \varrho_{n}M_{n}$. The last term is
typically small relative to the other two terms. Indeed, as noted
before, in many examples, for given $J_{n} \to\infty$ with
$J_{n}2^{J_{n}}/n = o(\tau_{J_{n}}^{2})$, $\epsilon_{n}$ can be taken
in such a way that $\epsilon_{n} \sim\sqrt{2^{J_{n}} (\log n)/n}$.
In that case $\tau_{J_{n}}^{-1} \epsilon_{n} \varrho_{n}M_{n} \sim
\tau_{J_{n}}^{-1}\varrho_{n} M_{n} \sqrt{2^{J_{n}}(\log n)/n} $, and
as long as $\varrho_{n} \to0$ sufficiently fast, that is, $\varrho
_{n} = o((\log n)^{-1/2})$, the convergence rate of the quasi-Bayes
estimator $\hat{g}_{\mathrm{QB}}$ reduces to $\max\{ 2^{-J_{n}s},\break \tau
_{J_{n}}^{-1}\sqrt{2^{J_{n}}/n} \}$.

\section{Prior specification: Examples}\label{sec4}
\label{secexamples}

In this section, we give some specific sieve priors for which the
quasi-posterior distribution (the quasi-Bayes estimator) attains
the minimax optimal rate of contraction (convergence, resp.). We
consider two types of priors, namely, product and isotropic priors.
We will verify that these priors meet conditions (P1) and (P2) in Theorem
\ref{thm1} with the choice (\ref{optimal-J}).
For the notational convenience, define
\[
\varepsilon_{n,s,r} = \cases{ n^{-s/(2s+2r+1)}, &\quad in the mildly
ill-posed case,
\vspace*{2pt}\cr
(\log n)^{-s}, &\quad in the severely ill-posed case.}
\]
We may think of the severely ill-posed case as the case with $r=\infty$.

\begin{proposition}
\label{prop1}
Suppose that Assumptions \ref{A1}--\ref{A4} are satisfied. Consider
the following two classes of prior distributions on $\R^{2^{J_{n}}}$:
\begin{longlist}[(Isotropic prior).]
\item[(Product prior).]
Let $q(x)$ be a probability density function on $\R$ such that for a
constant $A> \sup_{l \geq1} | b_{0l} |$:
(1) $q(x)$ is positive\vspace*{1pt} on $[-A,A]$;
(2) $\log q(x)$ is Lipschitz continuous on $[-A,A]$, that is, there
exists a constant $L > 0$ possibly depending on $A$ such that $| \log
q(x) - \log q(y) | \leq L |x-y|, \forall x,y \in[-A,A]$.
Take the\vspace*{1pt} density of the generating prior by $\tilde{\pi}_{n}
(b^{J_{n}}) = \prod_{l=1}^{2^{J_{n}}} q(b_{l})$.
\item[(Isotropic prior).] Let $r(x)$ be a probability density function
on $[0,\infty)$ having all moments such that: (1) for a constant $A> \|
g_{0} \|$, $r(x)$ is positive and continuous on $[0,A]$; (2) for a
constant $c'' > 0$, $\int_{0}^{\infty} x^{k-1} r(x) \,dx \leq e^{c'' k
\log k}$ for all $k$ sufficiently large.
Take the density\vspace*{-2pt} of the generating prior by $\tilde{\pi}_{n}
(b^{J_{n}}) \propto r(\| b^{J_{n}} \|_{\ell^{2}})$.
\end{longlist}

Take $J_{n}$ as in (\ref{optimal-J}). Then, in either case of product
or isotropic priors, for every sequence $M_{n} \to\infty$, we have
$\Pi_{n} \{ g\dvtx  \| g - g_{0} \| > M_{n} \varepsilon_{n,s,r} \mid\D
_{n} \} \stackrel{P}{\to} 0$. Furthermore, if $\sup_{w \in[0,1]} \E
[ U^{2} 1(| U | > \lambda) \mid W=w] \to0$ as $\lambda\to\infty$,
then there exists a constant $M > 0$ such that
$\Pi_{n} \{ g\dvtx  \| g - g_{0} \| > M \varepsilon_{n,s,r} \mid\D_{n}
\} \stackrel{P}{\to} 0$.
\end{proposition}
\begin{pf}
See Appendix D 
in the supplemental file \cite{K13}.
\end{pf}

Proposition \ref{prop1} shows that a wide class of priors constructed
on slowly growing sieves lead to the minimax optimal contraction rate
(see Remark~\ref{minimax}).
In either case of product or isotropic priors, the constant $A$ is not
necessarily known, which allows $q(x)$ and $r(x)$ to have unbounded support.
For example, in the former case, $q(x)$ may be the density of the
standard normal distribution, in which case $A$ can be taken to be
arbitrarily large.
Likewise, in the latter case, $r(x)$ may be the density of an
exponential distribution: $r(x) = \lambda e^{-\lambda x}, x \geq0$ for
some $\lambda> 0$. In the isotropic prior case, $r(x)$ should have all
moments, that is, $\int_{0}^{\infty} x^{k} r(x) \,dx < \infty$ for all
$k \geq1$, which ensures that $\tilde{\pi}_{n}(b^{J_{n}}) \propto
r(\| b^{J_{n}} \|_{\ell^{2}})$ is a proper distribution on $\R
^{2^{J_{n}}}$ for every $n \geq1$.

The next proposition shows that two classes of priors in Proposition
\ref{prop1} lead to the minimax optimal convergence rate for
the quasi-Bayes estimator.

\begin{proposition}
\label{prop2}
Suppose that Assumptions \ref{A1}--\ref{A4} are satisfied.
Furthermore, assume that $\sup_{w \in[0,1]} \E[ U^{2} 1(| U | >
\lambda) \mid W=w] \to0$ as $\lambda\to\infty$. Consider the two
classes of prior distributions on $\R^{2^{J_{n}}}$ given in
Proposition \ref{prop1}. In the isotropic prior case, assume further
that $r(x)$ is Lipschitz continuous on $[0,A]$. Take $J_{n}$ as in
(\ref{optimal-J}). Then, in either case of product or isotropic
priors, there exists a constant $M > 0$ such that
$\Prob\{ \| \hat{g}_{\mathrm{QB}} - g_{0} \| > M \varepsilon_{n,s,r} \} \to0$.
\end{proposition}

\begin{pf}
See Appendix D 
in the supplemental file \cite{K13}.
\end{pf}

\begin{remark}
In the above propositions, $J_{n}$ plays the role of regularization and
should be chosen sufficiently slowly growing, thereby there is no need
to place restrictions on weights on $b_{l}$ between $1 \leq l \leq2^{J_{n}}$.
The abstract Theorem \ref{thm1} is derived to cover this case. There
is another way to deal with the ill-posedness, that is, allowing for
large-dimensional sieves but placing prior distributions that have
smaller weights on $b_{l}$ for larger $l$ (``shrinking priors''), which
corresponds to the ``sieve method using large-dimensional sieves with
heavy penalties'' in the classification of \cite{CP11}.\footnote{The
previous version of this paper contains results on shrinking priors,
but $J_{n}$ should be still slowly growing as in the above
propositions, which corresponds to the sieve method using slowly
growing sieves with \emph{light} penalties. Those results have been removed in
the current version according to the referee's suggestion, but
available upon request.} The supplementary material of \cite{LJ11} is
concerned with this approach, but they did not establish sharp
contraction rates.
The extension to this approach requires a different technique than that
used in the present paper, and remains as an open problem.
\end{remark}

\section{\texorpdfstring{Proofs of Theorems \protect\ref{thm1} and \protect\ref{thm3}}
{Proofs of Theorems 1 and 2}}\label{sec5}

\subsection{\texorpdfstring{Proof of Theorem \protect\ref{thm1}}{Proof of Theorem 1}}\label{sec5.1}
\label{secprf1}

Before proving Theorem \ref{thm1}, we first prepare some technical
lemmas (Lemmas \ref{lemB1}--\ref{lemB4}) and establish preliminary
rates of contraction for the quasi-posterior distribution (Proposition
\ref{propA1}).
Proofs of Lemmas \ref{lemB1}--\ref{lemB4} are given in Appendix B %
in the supplemental file \cite{K13}.
For the notational convenience, define the matrices
\begin{eqnarray*}
\hat{\Phi}_{WX} &=& \En\bigl[ \phi^{J_{n}}(W_{i}) \phi
^{J_{n}}(X_{i})^{T}\bigr],\qquad \hat{\Phi}_{XW}
= \hat{\Phi}_{WX}^{T},\\
\hat{\Phi}_{WW}&=&\En\bigl[
\phi^{J_{n}}(W_{i})^{\otimes2}\bigr],
\end{eqnarray*}
which are the empirical counterparts of $\Phi^{WX},\Phi^{XW}$ and
$\Phi_{WW}$, respectively.
Also define
\[
U_{i}=Y_{i}-g_{0}(X_{i}),\qquad
R_{i}=Y_{i}-P_{J_{n}}g_{0}(X_{i}),\qquad
\Delta_{n}=\sqrt{n}\En\bigl[ \phi^{J_{n}}(W_{i})R_{i}
\bigr].
\]
Lemma \ref{lemB1} is a technical lemma on these quantities. Lemma \ref
{lemB2} characterizes the total variation convergence between two
centered multivariate normal distributions with increasing dimensions
in terms of the speed of convergence between the corresponding
covariance matrices.
Lemma \ref{lemB4} will be used in the latter part in the proof of
Theorem \ref{thm1}.
%
\begin{lemma}
\label{lemB1}
Suppose that Assumptions \ref{A1}--\ref{A4} are satisfied. Let $J_{n}
\to\infty$ as $n \to\infty$.
\textup{(i)} There exists a constant $D > 0$ such that $\sup_{w \in[0,1]} \|
\phi^{J} (w) \|_{\ell^{2}} \leq D 2^{J/2}$ for all $J \geq J_{0}$.
\textup{(ii)} $C_{1}^{-1} \leq s_{\min} (\E[ \phi^{J} (W)^{\otimes2} ] )
\leq s_{\max}(\E[ \phi^{J} (W)^{\otimes2} ]) \leq C_{1}$ and
$s_{\max} (\E[ \phi^{J}(W)\phi^{J}(X)^{T}]) \leq C_{1}$ for all $J
\geq J_{0}$. \textup{(iii)} If $J_{n} 2^{J_{n}}/n \to0$, $\| \hat{\Phi}_{WW}
- \Phi_{WW} \|_{\op} =O_{P}(\sqrt{J_{n}2^{J_{n}}/n})$ and $\| \hat
{\Phi}_{WX} - \Phi_{WX} \|_{\op} =O_{P}(\sqrt{J_{n}2^{J_{n}}/n})$.
\textup{(iv)} $\| \En[ \phi^{J_{n}}(W_{i})R_{i}] \|_{\ell^{2}}^{2} =
O_{P}(2^{J_{n}}/n + \tau_{J_{n}}^{2}2^{-2J_{n}s})$. \textup{(v)} If $J_{n}
2^{J_{n}}/n = o(\tau_{J_{n}}^{2})$, $s_{\min}(\hat{\Phi}_{WX}) \geq
(1-o_{P}(1))\tau_{J_{n}}$.
\end{lemma}

\begin{lemma}
\label{lemB2}
Let $\Sigma_{n}$ be a sequence of symmetric positive definite matrices
of dimension $k_{n} \to\infty$ as $n \to\infty$ such that $\|
\Sigma_{n} - I_{k_{n}} \|_{\op}=o(k_{n}^{-1})$.
Then as $n \to\infty$,
\[
\int\bigl| dN(0,\Sigma_{n}) (x) - dN(0,I_{k_{n}}) (x) \bigr| \,dx \to0.
\]
\end{lemma}

\begin{lemma}
\label{lemB4}
Let $\hat{A}_{n}$ be a sequence of random $k_{n} \times k_{n}$
matrices where $k_{n}$ is either bounded or $k_{n} \to\infty$ as $n
\to\infty$. Suppose that there exist sequences of positive constants
$\epsilon_{n}, \delta_{n}$ and a sequence of nonrandom,\vspace*{1pt} nonsingular
$k_{n} \times k_{n}$ matrices $A_{n}$ such that $\epsilon_{n} \to0,
\delta_{n} \to0, s_{\min}(A_{n}) \gtrsim\epsilon_{n}, \| \hat
{A}_{n} - A_{n} \|_{\op} = O_{P}(\delta_{n})$ and $\epsilon_{n}^{-1}
\delta_{n} \to0$. Then $\hat{A}_{n}$ is nonsingular with
probability approaching one and $\| \hat{A}_{n}^{-1} A_{n} - I_{k_{n}}
\|_{\op} \vee\| A_{n} \hat{A}_{n}^{-1} - I_{k_{n}} \|_{\op} =
O_{P}(\epsilon_{n}^{-1} \delta_{n})$.
\end{lemma}

The following proposition gives preliminary rates of contraction for
the quasi-posterior distribution.

\begin{proposition}[(Preliminary contraction rates)]
\label{propA1}
Suppose that Assumptions \ref{A1}--\ref{A4} are satisfied.
Take $J_{n}$ in such a way that $J_{n} \to\infty$ and $J_{n}
2^{J_{n}}/n = o(\tau_{J_{n}}^{2})$. Let $\epsilon_{n}$ be a sequence
of positive constants such that $\epsilon_{n} \to0$ and $\sqrt{n}
\epsilon_{n} \to\infty$.
Assume that a sequence of generating priors $\tilde{\Pi}_{n}$
satisfies condition \textup{(P1)} of Theorem \ref{thm1}. Define the
data-dependent, empirical seminorm $\| \cdot\|_{\D_{n}}$ on $\R
^{2^{J_{n}}}$ by
\[
\bigl\| b^{J_{n}} \bigr\|_{\D_{n}} = \bigl\| \hat{\Phi}_{WX}
b^{J_{n}} \bigr\|_{\ell
^{2}},\qquad b^{J_{n}} \in\R^{2^{J_{n}}}.
\]
Then for every sequence $M_{n} \to\infty$, we have
\[
\tilde{\Pi}_{n} \bigl\{ b^{J_{n}}\dvtx  \bigl\| b^{J_{n}} -
b^{J_{n}}_{0} \bigr\|_{\D
_{n}} > M_{n} \bigl(
\epsilon_{n} + \tau_{J_{n}} 2^{-J_{n}s}\bigr) \mid
\D_{n} \bigr\} \stackrel{P} {\to} 0.
\]
\end{proposition}

\begin{pf}{Proof of Proposition \ref{propA1}}
The proof consists of constructing suitable ``tests'' and is
essentially similar to, for example, the proof of Theorem 2.1 in \cite{GGV00}.
Let $\delta_{n} = \epsilon_{n} + \tau_{J_{n}} 2^{-J_{n}s}$. We wish
to show that there exists a constant $c_{0}>0$ such that
%
\begin{equation}\label{step2}
\Prob\bigl\{ \tilde{\Pi}_{n} \bigl( b^{J_{n}}\dvtx  \bigl\|
b^{J_{n}} - b_{0}^{J_{n}} \bigr\|_{\D_{n}} >
M_{n} \delta_{n} \mid\D_{n}\bigr) \leq
e^{-c_{0}M_{n}^{2}n\delta_{n}^{2}} \bigr\} \to1.
\end{equation}
Note that since $\sqrt{n} \epsilon_{n} \to\infty$, $n \delta
^{2}_{n} \geq n \epsilon_{n}^{2} \to\infty$.
Below, $c_{1},c_{2},\ldots$ are some positive constants of which the
values are understood in the context.\vspace*{1pt}

Note that $Y_{i}=P_{J_{n}} g_{0}(X_{i})+ R_{i} = \phi
^{J_{n}}(X_{i})^{T}b_{0}^{J_{n}} + R_{i}$. Then for $b^{J_{n}} \in\R
^{2^{J_{n}}}$,
%
\begin{eqnarray}\label{moment2}
\En\bigl[\hat{m}^{2}\bigl(W_{i},b^{J_{n}}\bigr)
\bigr] &=& -2\bigl(b^{J_{n}}-b^{J_{n}}_{0}
\bigr)^{T} \hat{\Phi}_{XW} \hat{\Phi}_{WW}^{-1}
\En\bigl[ \phi^{J_{n}}(W_{i})R_{i}\bigr]
\nonumber
\\
&&{} +\bigl(b^{J_{n}}-b^{J_{n}}_{0}\bigr)^{T}
\hat{\Phi}_{XW} \hat{\Phi}_{WW}^{-1} \hat{
\Phi}_{WX} \bigl(b^{J_{n}}-b^{J_{n}}_{0}\bigr)
\\
&&{} +\En\bigl[ \phi^{J_{n}}(W_{i})R_{i}
\bigr]^{T} \hat{\Phi}_{WW}^{-1} \En\bigl[
\phi^{J_{n}}(W_{i})R_{i}\bigr].\nonumber
\end{eqnarray}
Since the last term is independent of $b^{J_{n}}$, it is canceled out
in the quasi-posterior distribution. Denote by $\ell_{b^{J_{n}}}(\D
_{n})$ the sum of the first two terms in (\ref{moment2}). Then
\[
\tilde{\Pi}_{n}\bigl(db^{J_{n}} \mid\D_{n}\bigr)
\propto\exp\bigl\{ -(n/2) \ell_{b^{J_{n}}}(\D_{n}) \bigr\} \tilde{
\Pi}_{n}\bigl(db^{J_{n}}\bigr).
\]
Using the fact that for any $x,y,c \in\R$ with $c > 0$, $2xy \leq c
x^{2} + c^{-1} y^{2}$, we have
%
\begin{eqnarray}\label{lower}
\ell_{b^{J_{n}}}(\D_{n}) &\geq& (\hat{\lambda}_{\min} - c)\bigl\|
b^{J_{n}}-b^{J_{n}}_{0} \bigr\|_{\D_{n}}^{2}
\nonumber\\[-8pt]\\[-8pt]
&&{} - c^{-1}\hat{\lambda}_{\max}^{2} \bigl\| \En\bigl[ \phi
^{J_{n}}(W_{i})R_{i}\bigr] \bigr\|_{\ell^{2}}^{2}\qquad
\forall c > 0,\nonumber
\end{eqnarray}
where $\hat{\lambda}_{\min}$ and $\hat{\lambda}_{\max}$ are the
minimum and maximum eigenvalues of the matrix $\hat{\Phi}_{WW}^{-1}$,
respectively. Likewise, we have
%
\begin{eqnarray}\label{upper}
\ell_{b^{J_{n}}}(\D_{n}) &\leq& (\hat{\lambda}_{\max} + c)\bigl\|
b^{J_{n}}-b^{J_{n}}_{0} \bigr\|_{\D_{n}}^{2}
\nonumber\\[-8pt]\\[-8pt]
&&{} + c^{-1}\hat{\lambda}_{\max}^{2} \bigl\| \En\bigl[ \phi
^{J_{n}}(W_{i})R_{i}\bigr] \bigr\|_{\ell^{2}}^{2}\qquad
\forall c > 0.\nonumber
\end{eqnarray}

Define the event
\begin{eqnarray*}
\mathcal{E}_{1n} &=& \bigl\{ \D_{n}\dvtx  \hat{
\lambda}_{\min} < 0.5C_{1}^{-1} \bigr\} \cup\{
\D_{n}\dvtx  \hat{\lambda}_{\max} > 1.5 C_{1} \}
\\
&&{}\cup\bigl\{ \D_{n}\dvtx  \bigl\| \En\bigl[ \phi^{J_{n}}(W_{i})R_{i}
\bigr] \bigr\|_{\ell
^{2}}^{2} > M_{n} \delta_{n}^{2}
\bigr\}.
\end{eqnarray*}
Construct the ``tests'' $\omega_{n}$ by $\omega_{n} = 1(\mathcal
{E}_{1n})$. Then we have
%
\begin{eqnarray}\label{bound1}
&&\tilde{\Pi}_{n}\bigl(b^{J_{n}}\dvtx  \bigl\| b^{J_{n}} -
b_{0}^{J_{n}} \bigr\|_{\D
_{n}} > M_{n}
\delta_{n} \mid\D_{n}\bigr)
\nonumber
\\
&&\qquad= \tilde{\Pi}_{n}\bigl(b^{J_{n}}\dvtx  \bigl\| b^{J_{n}} -
b_{0}^{J_{n}} \bigr\|_{\D
_{n}} > M_{n}
\delta_{n} \mid\D_{n}\bigr) \bigl\{ \omega_{n} + (1-
\omega_{n}) \bigr\}
\\
&&\qquad\leq\omega_{n} + \tilde{\Pi}_{n}\bigl(b^{J_{n}}\dvtx  \bigl\|
b^{J_{n}} - b_{0}^{J_{n}} \bigr\|_{\D_{n}} >
M_{n} \delta_{n} \mid\D_{n}\bigr) (1-\omega
_{n}).\nonumber
\end{eqnarray}
By Lemma \ref{lemB1}(ii)--(iv), we have $\Prob( \omega_{n}=1 ) =
\Prob(\mathcal{E}_{1n}) \to0$.

For the second term in (\ref{bound1}), taking $c > 0$ sufficiently
small in (\ref{lower}), we have
\begin{eqnarray*}
&&(1-\omega_{n}) \int_{\| b^{J_{n}} - b_{0}^{J_{n}} \|_{\D_{n}} >
M_{n} \delta_{n}} \exp\bigl\{ -(n/2)
\ell_{b^{J_{n}}}(\D_{n}) \bigr\} \tilde{\Pi}_{n}
\bigl(db^{J_{n}}\bigr)
\\
&&\qquad\leq\exp\bigl\{ -c_{1}M_{n}^{2}n
\delta_{n}^{2} + O\bigl( M_{n}n\delta
_{n}^{2} \bigr) \bigr\} \leq e^{-c_{2} M_{n}^{2} n\delta_{n}^{2}}.
\end{eqnarray*}
On the other hand, taking, say $c=1$ in (\ref{upper}), we have
\begin{eqnarray*}
&&(1-\omega_{n}) \int\exp\bigl\{ -(n/2) \ell_{b^{J_{n}}}(
\D_{n}) \bigr\} \tilde{\Pi}_{n}\bigl(db^{J_{n}}\bigr)
\\
&&\qquad\geq(1-\omega_{n}) \int_{\| b^{J_{n}} - b_{0}^{J_{n}} \|_{\D_{n}}
\leq\sqrt{M_{n}} \epsilon_{n}} \exp\bigl\{ -(n/2)
\ell_{b^{J_{n}}}(\D_{n}) \bigr\} \tilde{\Pi}_{n}
\bigl(db^{J_{n}}\bigr)
\\
&&\qquad\geq(1-\omega_{n}) e^{-c_{3}M_{n}n\delta^{2}_{n}} \int_{\|
b^{J_{n}} - b_{0}^{J_{n}} \|_{\D_{n}} \leq\sqrt{M_{n}} \epsilon
_{n}}
\tilde{\Pi}_{n}\bigl(db^{J_{n}}\bigr).
\end{eqnarray*}
Denote by $\hat{s}_{\max}$ the maximum singular value of the matrix
$\hat{\Phi}_{WX}$, so that
\[
\bigl\| b^{J_{n}} - b_{0}^{J_{n}} \bigr\|_{\D_{n}} \leq
\hat{s}_{\max} \bigl\| b^{J_{n}} - b_{0}^{J_{n}}
\bigr\|_{\ell^{2}}.
\]
Define the event $\mathcal{E}_{2n} = \{ \D_{n}\dvtx  \hat{s}_{\max}
\leq1.5 C_{1} \}$. By Lemma \ref{lemB1}(ii) and (iii), we have
$\Prob(\mathcal{E}_{2n}) \to1$.
Since $M_{n} \to\infty$, for all $n$ sufficiently large, we have
\begin{eqnarray*}
&&1(\mathcal{E}_{2n}) (1-\omega_{n}) \int\exp\bigl\{ -(n/2)
\ell_{b^{J_{n}}}(\D_{n}) \bigr\} \tilde{\Pi}_{n}
\bigl(db^{J_{n}}\bigr)
\\
&&\qquad\geq1(\mathcal{E}_{2n}) (1-\omega_{n}) e^{-c_{3}M_{n}n\delta
^{2}_{n}}
\tilde{\Pi}_{n}\bigl(b^{J_{n}}\dvtx  \bigl\| b^{J_{n}} -
b_{0}^{J_{n}} \bigr\| _{\ell{2}} \leq\epsilon_{n}
\bigr)
\\
&&\qquad\geq1(\mathcal{E}_{2n}) (1-\omega_{n}) e^{-c_{3}M_{n}\delta
^{2}_{n} - Cn\epsilon_{n}^{2}}
\\
&&\qquad\geq1(\mathcal{E}_{2n}) (1-\omega_{n}) e^{-c_{4}M_{n}n\delta^{2}_{n}},
\end{eqnarray*}
where the second inequality is due to the small ball condition (P1).
Summarizing, we have
\[
\tilde{\Pi}_{n}\bigl(b^{J_{n}}\dvtx  \bigl\| b^{J_{n}} -
b_{0}^{J_{n}} \bigr\|_{\D_{n}} > M_{n}
\delta_{n} \mid\D_{n}\bigr) (1-\omega_{n}) \leq1
\bigl(\mathcal{E}_{2n}^{c}\bigr) + e^{-c_{2}M_{n}^{2} n \delta_{n}^{2}
+ c_{4} M_{n} n
\delta^{2}_{n}}.
\]
Therefore, we obtain (\ref{step2}) for a sufficiently small $c_{0}>0$.
\end{pf}

We are now in position to prove Theorem \ref{thm1}. We will say that a
sequence of random variables $A_{n}$ is \emph{eventually} bounded by
another sequence of random variables $B_{n}$ if $\Prob(A_{n} \leq
B_{n}) \to1$ as $n \to\infty$.

\begin{pf*}{Proof of Theorem \ref{thm1}}
We first note that by Lemma \ref{lemB1}(ii), (iii) and (v), the
matrices $\hat{\Phi}_{WX}$ and $\hat{\Phi}_{WW}$ are nonsingular with
probability approaching one. Conditional on $\D_{n}$, define the
rescaled ``parameter'' $\theta
^{J_{n}}=(\theta_{1},\ldots,\theta_{2^{J_{n}}})^{T} = \sqrt{n}\hat
{\Phi}_{WX}(b^{J_{n}} - b_{0}^{J_{n}})$. By (\ref{moment2}), the
corresponding ``quasi-posterior'' density for $\theta^{J_{n}}$ is given
by
\[
\pi^{*}_{n}\bigl(\theta^{J_{n}} \mid\D_{n}
\bigr)\,d\theta^{J_{n}} \propto\tilde{\pi}_{n}
\bigl(b_{0}^{J_{n}} + \hat{\Phi}_{WX}^{-1}
\theta^{J_{n}}/\sqrt{n}\bigr) dN(\Delta_{n},\hat{
\Phi}_{WW}) \bigl(\theta^{J_{n}}\bigr) \,d\theta^{J_{n}},
\]
where recall that $\Delta_{n}=\sqrt{n} \En[ \phi^{J_{n}} (W_{i})
R_{i} ]$ (this operation is valid as soon as $\hat{\Phi}_{WX}$ and
$\hat{\Phi}_{WW}$ are nonsingular, of which the probability is
approaching one).

The proof of Theorem \ref{thm1} consists of 3 steps. After step 1, we
will turn to the proof of (\ref{rate}).
The remaining two steps are devoted to the proof of (\ref{BvM}).

Step 1. We first show that
%
\begin{equation}\label{2step1}
\int\bigl| \pi^{*}_{n}\bigl(\theta^{J_{n}} \mid
\D_{n}\bigr) - dN(\Delta_{n},\hat{\Phi}_{WW})
\bigl(\theta^{J_{n}}\bigr) \bigr| \,d\theta^{J_{n}} \stackrel{P} {\to} 0.
\end{equation}
In this step, we do \emph{not} assume $J_{n}2^{3J_{n}}/n = o(\tau_{J_{n}}^{2})$.
As before, let $\delta_{n} = \epsilon_{n} + \tau_{J_{n}}
2^{-J_{n}s}$. By Proposition \ref{propA1}, for every sequence $M_{n}
\to\infty$,
\[
\int_{\| \theta^{J_{n}} \|_{\ell^{2}} \leq M_{n} \sqrt{n} \delta
_{n}} \pi^{*}_{n}\bigl(
\theta^{J_{n}} \mid\D_{n}\bigr)\,d\theta^{J_{n}} = 1 +
o_{P}(1)
\]
by which we have
%
\begin{eqnarray}\label{bound2}
&&\mbox{Left-hand side of (\ref{2step1})}
\nonumber
\\
&&\qquad\leq\int_{\| \theta^{J_{n}} \|_{\ell^{2}} \leq M_{n} \sqrt{n}
\delta_{n}} \bigl| \pi^{*}_{n}\bigl(
\theta^{J_{n}} \mid\D_{n}\bigr) - dN(\Delta_{n},\hat{
\Phi}_{WW}) \bigl(\theta^{J_{n}}\bigr) \bigr| \,d\theta^{J_{n}}
\\
&&\qquad\quad{} + \int_{\| \theta^{J_{n}} \|_{\ell^{2}} > M_{n} \sqrt{n}
\delta_{n}}dN(\Delta_{n},\hat{
\Phi}_{WW}) \bigl(\theta^{J_{n}}\bigr) \,d\theta^{J_{n}} +
o_{P}(1).\nonumber
\end{eqnarray}
By Lemma \ref{lemB1}(iv), $\| \Delta_{n} \|_{\ell^{2}}=O_{P}(\sqrt
{n}\delta_{n})$, and by Lemma \ref{lemB1}(ii) and (iii),
$(1-o_{P}(1))C_{1}^{-1} \leq s_{\min} (\hat{\Phi}_{WW}) \leq s_{\max
}(\hat{\Phi}_{WW}) \leq(1+o_{P}(1))C_{1}$, so that the second
integral is eventually bounded by
%
\begin{equation}\label{integral}
\int_{\| \theta^{J_{n}} \|_{\ell^{2}} > \sqrt{M_{n} n} \delta
_{n}}dN(0,I_{2^{J_{n}}}) \bigl(\theta^{J_{n}}
\bigr) \,d\theta^{J_{n}},
\end{equation}
where note that $M_{n}$ is replaced by $\sqrt{M_{n}}$ to ``absorb''
the constant. By Borell's inequality for Gaussian measures (see, e.g.,
\cite{VW96}, Lemma A.2.2), for every $x > 0$,
%
\begin{equation}\label{borell}
\Prob\bigl( \bigl\| N(0,I_{2^{J_{n}}}) \bigr\|_{\ell^{2}} > \sqrt{2^{J_{n}}}
+ x \bigr) \leq e^{-x^{2}/2}.
\end{equation}
Here since $n\delta_{n}^{2} \geq n \epsilon_{n}^{2} \gtrsim
2^{J_{n}}$, $\sqrt{M_{n}n} \delta_{n}/\sqrt{2^{J_{n}}} \to\infty$,
so that the integral in (\ref{integral}) is $o(1)$.

It remains to show that the first integral in (\ref{bound2}) is
$o_{P}(1)$. This step uses a standard cancellation argument. Let $\C
_{n}:= \{ \theta^{J_{n}} \in\R^{2^{J_{n}}}\dvtx
\| \theta^{J_{n}} \|_{\ell^{2}} \leq M_{n} \sqrt{n} \delta_{n} \}
$. First, provided that $\| \hat{\Phi}^{-1}_{WX} \|_{\op} \leq1.5
\tau_{J_{n}}^{-1}$, for all $\theta^{J_{n}} \in\C_{n}$,
\[
\bigl\| \hat{\Phi}_{WX}^{-1}\theta^{J_{n}}/\sqrt{n}
\bigr\|_{\ell^{2}} \leq1.5 M_{n} \tau_{J_{n}}^{-1}
\delta_{n} \leq1.5 M_{n} \bigl(2^{-J_{n}s} +
\tau_{J_{n}}^{-1} \epsilon_{n}\bigr) \sim
M_{n} \gamma_{n}.
\]
So taking $M_{n} \to\infty$ such that $M_{n} = o(L_{n})$, $\| \hat
{\Phi}_{WX}^{-1}\theta^{J_{n}}/\sqrt{n} \|_{\ell^{2}} \leq L_{n}
\gamma_{n}$ and hence $\tilde{\pi}_{n}(b_{0}^{J_{n}} + \hat{\Phi
}_{WX}^{-1}\theta^{J_{n}}/\sqrt{n}) > 0$ for all $n$ sufficiently large.
Here, by Lemma \ref{lemB1}(v), we have $\Prob(\| \hat{\Phi
}^{-1}_{WX} \|_{\op} \leq1.5 \tau_{J_{n}}^{-1}) \to1$.

Suppose that $\| \hat{\Phi}^{-1}_{WX} \|_{\op} \leq1.5 \tau
_{J_{n}}^{-1}$. Let
\[
\pi^{*}_{n,\C_{n}}\bigl(\theta^{J_{n}} \mid
\D_{n}\bigr) \quad\mbox{and}\quad dN^{\C_{n}}(\Delta_{n},\hat{
\Phi}_{WW}) \bigl(\theta^{J_{n}}\bigr)
\]
denote the probability densities obtained by first restricting
$\pi^{*}_{n}(\theta^{J_{n}} \mid\D_{n})$ and $dN(\Delta_{n},\hat
{\Phi}_{WW})(\theta^{J_{n}})$ to the ball $\C_{n}$ and then
renormalizing, respectively. By the first part of the present proof,
replacing $\pi^{*}_{n}(\theta^{J_{n}} \mid\D_{n})$ and $dN(\Delta
_{n},\hat{\Phi}_{WW})(\theta^{J_{n}})$ by $\pi^{*}_{n,\C
_{n}}(\theta^{J_{n}} \mid\D_{n})$ and $dN^{\C_{n}}(\Delta_{n},\hat
{\Phi}_{WW})(\theta^{J_{n}})$, respectively, in the first integral in
(\ref{bound2}) has impact at most $o_{P}(1)$.
Abbreviating $\pi^{*}_{n,\C_{n}}(\theta^{J_{n}} \mid\D_{n})$ by
$\pi^{*}_{n,\C_{n}}$, $dN^{\C_{n}}(\Delta_{n},\hat{\Phi
}_{WW})(\theta^{J_{n}})$ by $dN^{\C_{n}}$, $dN(\Delta_{n},\hat{\Phi
}_{WW})(\theta^{J_{n}})$ by $dN$, and $\tilde{\pi}_{n}(b_{0}^{J_{n}}
+ \hat{\Phi}_{WX}^{-1}\theta^{J_{n}}/\sqrt{n})$ by $\tilde{\pi
}_{n}$, we have
\begin{eqnarray*}
\int\bigl| \pi^{*}_{n,\C_{n}} - dN^{\C_{n}}\bigr| &=& \int\biggl
\llvert1 - \frac
{dN^{\C_{n}}}{\pi^{*}_{n,\C_{n}}} \biggr\rrvert\pi^{*}_{n,\C_{n}} =
\int\biggl\llvert1 - \frac{dN/\int_{\C_{n}} dN}{ \tilde{\pi}_{n} dN
/\int_{\C_{n}} \tilde{\pi}_{n} dN } \biggr\rrvert\pi^{*}_{n,\C_{n}}
\\
&=&\int\biggl\llvert1 - \frac{\int_{\C_{n}}\tilde{\pi}_{n} dN }{ \tilde
{\pi}_{n}\int_{\C_{n}} dN} \biggr\rrvert\pi^{*}_{n,\C_{n}}
=\int\biggl\llvert1 - \frac{\int_{\C_{n}}\tilde{\pi}_{n} dN^{\C_{n}}}{
\tilde
{\pi}_{n}} \biggr\rrvert\pi^{*}_{n,\C_{n}}.
\end{eqnarray*}
By the convexity of the map $x \mapsto|1-x|$ and Jensen's inequality,
the last expression is bounded by
\[
\sup_{\theta^{J_{n}} \in\C_{n}, \tilde{\theta}^{J_{n}} \in\C
_{n}} \biggl\llvert1 - \frac{\tilde{\pi}_{n}(b_{0}^{J_{n}} + \hat{\Phi
}_{WX}^{-1}\theta^{J_{n}}/\sqrt{n})}{\tilde{\pi}_{n}(b_{0}^{J_{n}}
+ \hat{\Phi}_{WX}^{-1}\tilde{\theta}^{J_{n}}/\sqrt{n})} \biggr\rrvert,
\]
which is eventually bounded by
\[
\sup_{\| b^{J_{n}} \|_{\ell^{2}} \leq L_{n} \gamma_{n}, \| \tilde
{b}^{J_{n}} \|_{\ell^{2}}\leq L_{n}\gamma_{n} } \biggl\llvert1 - \frac
{\tilde{\pi}_{n}(b_{0}^{J_{n}} + b^{J_{n}})}{\tilde{\pi
}_{n}(b_{0}^{J_{n}} +\tilde{b}^{J_{n}})} \biggr\rrvert.
\]
The last expression goes to zeros as $n \to\infty$ by condition (P2).

We now turn to the proof of (\ref{rate}).
Take any $M_{n} \to\infty$ (this $M_{n}$ may be different from the
previous $M_{n}$). By step 1, we have
\begin{eqnarray*}
&&
\sup_{z > 0}  \biggl| \tilde{\Pi}_{n}\bigl\{
b^{J_{n}}\dvtx  \bigl\| \hat{\Phi}_{WX}\bigl(b^{J_{n}}-b_{0}^{J_{n}}
\bigr) \bigr\|_{\ell^{2}} > z \mid\D_{n} \bigr\}
\\
&&\qquad{} - \int_{\| \theta^{J_{n}} \|_{\ell^{2}} > z} dN\bigl(n^{-1/2}\Delta_{n},n^{-1}
\hat{\Phi}_{WW}\bigr) \bigl(\theta^{J_{n}}\bigr) \,d
\theta^{J_{n}} \biggr| \stackrel{P} {\to} 0.
\end{eqnarray*}
By Lemma \ref{lemB1}(v), we have
\begin{eqnarray*}
\bigl\| \hat{\Phi}_{WX}\bigl(b^{J_{n}}-b_{0}^{J_{n}}
\bigr) \bigr\|_{\ell^{2}} &\geq& s_{\min}(\hat{\Phi}_{WX}) \bigl\|
b^{J_{n}}-b_{0}^{J_{n}} \bigr\|_{\ell^{2}}
\\
&\geq& \bigl(1-o_{P}(1)\bigr) \tau_{J_{n}} \bigl\|
b^{J_{n}}-b_{0}^{J_{n}} \bigr\|_{\ell^{2}}
\end{eqnarray*}
by which we have, uniformly in $z > 0$,
\begin{eqnarray*}
&&\tilde{\Pi}_{n}\bigl\{ b^{J_{n}}\dvtx  \bigl\| b^{J_{n}}-b_{0}^{J_{n}}
\bigr\|_{\ell
^{2}} > 2 \tau_{J_{n}}^{-1} z \mid\D_{n}
\bigr\}
\\
&&\qquad\leq\tilde{\Pi}_{n}\bigl\{ b^{J_{n}}\dvtx  \bigl\| \hat{\Phi
}_{WX}\bigl(b^{J_{n}}-b_{0}^{J_{n}}\bigr)
\bigr\|_{\ell^{2}} > z \mid\D_{n} \bigr\} + o_{P}(1)
\\
&&\qquad\leq\int_{\| \theta^{J_{n}} \|_{\ell^{2}} > z} dN\bigl(n^{-1/2}\Delta
_{n},n^{-1}\hat{\Phi}_{WW}\bigr) \bigl(
\theta^{J_{n}}\bigr) \,d\theta^{J_{n}} + o_{P}(1).
\end{eqnarray*}
By Markov's inequality, the integral in the last expression is bounded by
\[
\frac{1}{nz^{2}} \bigl\{ \| \Delta_{n} \|_{\ell^{2}}^{2}
+ \operatorname{tr} (\hat{\Phi}_{WW}) \bigr\}.
\]
By Lemma \ref{lemB1}(ii)--(iv), we have $ \| \Delta_{n} \|_{\ell
^{2}}^{2} + \operatorname{tr} (\hat{\Phi}_{WW}) = O_{P}(2^{J_{n}} + n\tau
_{J_{n}}^{2} 2^{-2J_{n}s})$. Therefore, we conclude that, taking
$z=M_{n} (\tau_{J_{n}}2^{-J_{n}s} + \sqrt{2^{J_{n}}/n})$, $\tilde
{\Pi}_{n}\{ b^{J_{n}}\dvtx\break  \| b^{J_{n}}-b_{0}^{J_{n}} \|_{\ell^{2}} > 2
M_{n} (2^{-J_{n}s} + \tau_{J_{n}}^{-1} \sqrt{2^{J_{n}}/n}) \mid\D
_{n} \} \stackrel{P}{\to} 0$, which leads to
the contraction rate result (\ref{rate}).

In what follows, we assume $J_{n}2^{3J_{n}}/n = o(\tau_{J_{n}}^{2})$,
and prove the asymptotic normality result (\ref{BvM}).

Step 2 (replacement of $\hat{\Phi}_{WW}$ by $\Phi_{WW}$). This step
shows that
\[
\int\bigl| dN(\Delta_{n},\hat{\Phi}_{WW}) \bigl(
\theta^{J_{n}}\bigr) - dN(\Delta_{n},\Phi_{WW})
\bigl(\theta^{J_{n}}\bigr) \bigr| \,d\theta^{J_{n}} \stackrel{P} {\to} 0,
\]
which is equivalent to
\[
\int\bigl| dN(0,\hat{\Phi}_{WW}) \bigl(\theta^{J_{n}}\bigr) - dN(0,
\Phi_{WW}) \bigl(\theta^{J_{n}}\bigr) \bigr| \,d\theta^{J_{n}}
\stackrel{P} {\to} 0.
\]
By Lemmas \ref{lemB1}(ii), (iii) and \ref{lemB2}, this follows
if $\sqrt{J_{n}2^{J_{n}}/n} = o(2^{-J_{n}})$, that is, $J_{n} 2^{3J_{n}}
= o(n)$, which is satisfied since $J_{n}2^{3J_{n}}/n = o(\tau
_{J_{n}}^{2}) = o(1)$.\eject

Step 3 (replacement of $\hat{\Phi}_{WX}$ by $\Phi_{WX}$). We have
shown that
\[
\int\bigl| \pi^{*}_{n}\bigl(\theta^{J_{n}} \mid
\D_{n}\bigr) - dN(\Delta_{n},\Phi_{WW}) \bigl(
\theta^{J_{n}}\bigr) \bigr| \,d\theta^{J_{n}} \stackrel{P} {\to} 0.
\]
By Scheff\'{e}'s lemma, this means that
\[
\bigl\| \tilde{\Pi}_{n} \bigl\{ b^{J_{n}}\dvtx  \sqrt{n}\hat{\Phi
}_{WX}\bigl(b^{J_{n}} - b_{0}^{J_{n}}\bigr)
\in\cdot\mid\D_{n} \bigr\} - N(\Delta_{n},
\Phi_{WW}) (\cdot)\bigr\|_{\TV} \stackrel{P} {\to} 0
\]
or equivalently,
\[
\bigl\| \tilde{\Pi}_{n} \bigl\{ b^{J_{n}}\dvtx  \sqrt{n}
\bigl(b^{J_{n}} - b_{0}^{J_{n}}\bigr) \in\cdot\mid
\D_{n}\bigr\} - N\bigl(\hat{\Phi}_{WX}^{-1}
\Delta_{n},\hat{\Phi}_{WX}^{-1}\Phi_{WW}
\hat{\Phi}_{XW}^{-1}\bigr) (\cdot) \bigr\|_{\TV}
\stackrel{P} {\to} 0.
\]
The last expression is asymptotically valid since $\hat{\Phi}_{WX}$
is nonsingular with probability approaching one.
Recall the maximum quasi-likelihood estimator $\hat{b}^{J_{n}}$. With
probability approaching one, we have
\[
\hat{b}^{J_{n}} = \hat{\Phi}_{WX}^{-1} \En\bigl[
\phi^{J_{n}} (W_{i}) Y_{i} \bigr] =
b_{0}^{J_{n}} + \hat{\Phi}_{WX}^{-1} \En
\bigl[ \phi^{J_{n}} (W_{i}) R_{i} \bigr],
\]
so that $\sqrt{n} ( \hat{b}^{J_{n}} - b_{0}^{J_{n}}) = \hat{\Phi
}_{WX}^{-1} \Delta_{n}$.
Hence to conclude the theorem, it suffices to show that
%
\begin{eqnarray}\label{assertion1}
&&\bigl\| N\bigl(\hat{\Phi}_{WX}^{-1}\Delta_{n},\hat{
\Phi}_{WX}^{-1}\Phi_{WW}\hat{\Phi}_{XW}^{-1}
\bigr)
\nonumber\\[-8pt]\\[-8pt]
&&\qquad{} - N\bigl(\hat{\Phi}_{WX}^{-1}\Delta_{n},\Phi
_{WX}^{-1}\Phi_{WW}\Phi_{XW}^{-1}
\bigr) \bigr\|_{\TV} \stackrel{P} {\to} 0.\nonumber
\end{eqnarray}

Assertion (\ref{assertion1}) reduces to
\[
\bigl\| N\bigl(0,\Phi_{WX} \hat{\Phi}_{WX}^{-1}
\Phi_{WW}\hat{\Phi}_{XW}^{-1}\Phi_{XW}
\bigr) - N(0,\Phi_{WW}) \bigr\|_{\TV} \stackrel{P} {\to} 0.
\]
By Lemmas \ref{lemB1}(ii), (iii) and \ref{lemB4},
\[
\bigl\| \Phi
_{WX} \hat{\Phi}_{WX}^{-1}\Phi_{WW}\hat{\Phi}_{XW}^{-1}\Phi_{XW}
- \Phi_{WW} \bigr\|_{\op} = O_{P}\bigl(\tau_{J_{n}}^{-1} \sqrt
{J_{n}2^{J_{n}}/n}\bigr) = o_{P}\bigl(2^{-J_{n}}\bigr)
\]
[the last equality follows
since $J_{n}2^{3J_{n}}/n = o(\tau_{J_{n}}^{2})$]. Since $C_{1}^{-1}
\leq s_{\min}(\Phi_{WW}) \leq s_{\max}(\Phi_{WW}) \leq C_{1}$, the
desired conclusion follows from Lemma \ref{lemB2}.

Steps 1--3 lead to the asymptotic normality result (\ref{BvM}).
\end{pf*}

\subsection{\texorpdfstring{Proof of Theorem \protect\ref{thm3}}{Proof of Theorem 2}}\label{sec5.2}
\label{secprf3}

We first prove the following lemma.
%
\begin{lemma}
\label{lemprf3}
Suppose that the conditions of Theorem \ref{thm3} are satisfied. Then
there exists a constant $D > 0$ such that
\[
\Prob\bigl\{ \bigl\| \En\bigl[ \phi^{J_{n}}(W_{i})
U_{i} \bigr] \bigr\|_{\ell^{2}} > D \sqrt{2^{J_{n}}/n} \bigr\}
\to0.
\]
\end{lemma}

\begin{remark}
It is standard to show that $\| \En[ \phi^{J_{n}}(W_{i}) U_{i} ] \|
_{\ell^{2}} =\break O_{P}(\sqrt{2^{J_{n}}/n})$, which, however, does not
leads to the conclusion of Lemma \ref{lemprf3} since the former only
implies that for every sequence $M_{n} \to\infty$,\break $\Prob\{ \| \En[
\phi^{J_{n}}(W_{i}) U_{i} ] \|_{\ell^{2}} > M_{n} \sqrt{2^{J_{n}}/n}
\} \to0$. Hence, an additional step is needed. The current proof uses
a truncation argument and Talagrand's concentration inequality.
\end{remark}

\begin{pf*}{Proof of Lemma \ref{lemprf3}}
For a given $\lambda> 0$, define $U_{i}^{-} = U_{i} 1(|U_{i}| \leq
\lambda)$ and $U_{i}^{+} = U_{i} 1(| U_{i} | > \lambda)$.
Since $0 = \E[ U \mid W ] = \E[ U^{-} \mid W ] + \E[ U^{+} \mid W
]$, we have $\En[ \phi^{J_{n}}(W_{i}) U_{i} ] = n^{-1} \sum_{i=1}^{n}
\{ \phi^{J_{n}}(W_{i}) U^{-}_{i} - \E[ \phi^{J_{n}}(W)
U^{-}] \} +\break n^{-1} \*\sum_{i=1}^{n} \{ \phi^{J_{n}}(W_{i}) U^{+}_{i} -
\E[ \phi^{J_{n}}(W) U^{+}] \}$, by which we have
\begin{eqnarray*}
\bigl\| \En\bigl[ \phi^{J_{n}}(W_{i}) U_{i} \bigr]
\bigr\|_{\ell^{2}} &\leq& \Biggl\| n^{-1} \sum_{i=1}^{n}
\bigl\{ \phi^{J_{n}}(W_{i}) U^{-}_{i} -
\E\bigl[ \phi^{J_{n}}(W) U^{-}\bigr] \bigr\} \Biggr\|_{\ell^{2}}
\\
&&{} + \Biggl\| n^{-1}  \sum_{i=1}^{n} \bigl\{
\phi^{J_{n}}(W_{i}) U^{+}_{i} - \E\bigl[
\phi^{J_{n}}(W) U^{+}\bigr] \bigr\} \Biggr\|_{\ell^{2}}
\\
&=:& I + \mathit{II}.
\end{eqnarray*}

First, by Markov's inequality, we have for every $z > 0$,
\begin{eqnarray*}
\Prob( \mathit{II} > z ) &\leq& \frac{\E[ \mathit{II}^{2} ]}{z^{2}} \leq\frac{\sum
_{l=1}^{2^{J_{n}}} \E[ (\phi_{l}(W) U^{+})^{2}] }{n z^{2}}
\\
&\leq& \frac{\sup_{w \in[0,1]} \E[ U^{2} 1(| U | > \lambda) \mid
W=w] \times\sum_{l=1}^{2^{J_{n}}} \E[ \phi_{l}(W)^{2} ] }{n z^{2}}
\\
&\leq& \frac{C_{1}2^{J_{n}}}{n z^{2}} \times\sup_{w \in[0,1]} \E\bigl[
U^{2} 1\bigl(| U | > \lambda\bigr) \mid W=w\bigr],
\end{eqnarray*}
where we have used that $\sum_{l=1}^{2^{J_{n}}} \E[ \phi_{l}(W)^{2}
] = \tr(\Phi_{WW}) \leq2^{J_{n}} s_{\max}(\Phi_{WW}) \leq C_{1}
2^{J_{n}}$ by Lemma \ref{lemB1}(ii).
Thus, we have
\[
\Prob\bigl\{ \mathit{II} > \sqrt{C_{1}2^{J_{n}}/n} \bigr\} \leq\sup
_{w \in[0,1]} \E\bigl[ U^{2} 1\bigl(| U | > \lambda\bigr) \mid W=w
\bigr].
\]
By assumption, the right-hand side goes to zero as $\lambda\to\infty$.

Second, let $Z_{i} = \phi^{J_{n}}(W_{i}) U_{i}^{-} - \E[ \phi
^{J_{n}}(W) U^{-} ]$ (denote by $Z$ the generic version of $Z_{i}$).
Let $\mathbb{S}^{2^{J_{n}}-1}:= \{ \alpha^{J_{n}} \in\R
^{2^{J_{n}}}\dvtx  \| \alpha^{J_{n}} \|_{\ell^{2}} = 1 \}$. Then
\[
I = \bigl\| \En[ Z_{i} ] \bigr\|_{\ell^{2}} =\sup_{\alpha^{J_{n}} \in\mathbb
{S}^{2^{J_{n}}-1}}
\En\bigl[ \bigl(\alpha^{J_{n}}\bigr)^{T} Z_{i}\bigr].
\]
We make use of Talagrand's concentration inequality to bound the tail
probability of $I$. For any $\alpha^{J_{n}} \in\mathbb
{S}^{2^{J_{n}}-1}$, by Lemma \ref{lemB1}, we have
\begin{eqnarray*}
\E\bigl[ \bigl\{ \bigl(\alpha^{J_{n}}\bigr)^{T} Z \bigr
\}^{2} \bigr] &\leq&\sup_{w \in[0,1]} \E\bigl[ U^{2}
\mid W = w\bigr] \times s_{\max} (\Phi_{WW}) \leq
C_{1}^{2},
\\
\bigl|\bigl(\alpha^{J_{n}}\bigr)^{T} Z \bigr| &\leq&\lambda\sup
_{w \in[0,1]} \bigl\| \phi^{J_{n}}(w) \bigr\|_{\ell^{2}} \leq
D_{1} \lambda\sqrt{2^{J_{n}}}
\end{eqnarray*}
and
\begin{eqnarray*}
\bigl(\E[ I ]\bigr)^{2} &\leq&\E\bigl[ I^{2} \bigr] \leq
n^{-1} \sup_{w \in[0,1]} \E\bigl[ U^{2} \mid W = w
\bigr] \times\sum_{l=1}^{2^{J_{n}}} \E\bigl[
\phi_{l}(W)^{2} \bigr] \\
&\leq& C_{1}^{2}
2^{J_{n}}/n,
\end{eqnarray*}
where $D_{1} > 0$ is a constant.
Thus, by Talagrand's inequality (see Theorem 2 in
Appendix E
), we have for every $z > 0$
\[
\Prob\bigl\{ I \geq D_{2}\bigl( \sqrt{2^{J_{n}}/n} + \sqrt{z/n}
+ z \lambda\sqrt{2^{J_{n}}}/n \bigr) \bigr\} \leq e^{-z},
\]
where $D_{2} > 0$ is a constant independent of $\lambda$ and $z$.

The final conclusion follows from taking $\lambda= \lambda_{n} \to
\infty$ and $z=z_{n} \to\infty$ sufficiently slowly.
\end{pf*}

\begin{pf*}{Proof of Theorem \ref{thm3}}
Let $D_{1}$ and $D_{2}$ be some positive constants of which the values
are understood in the context.
For either $g_{0} \in B^{s}_{\infty,\infty}$ or $B^{s}_{2,2}$, $\|
g_{0} - P_{J_{n}} g_{0} \| = O(2^{-J_{n}s}) = o(1)$, by which we have
\begin{eqnarray*}
&&
\sum_{l=1}^{2^{J_{n}}} \Var\bigl\{ \En\bigl[
\phi_{l}(W_{i}) (g_{0}-P_{J_{n}}g_{0})
(X_{i}) \bigr] \bigr\} \\
&&\qquad\leq n^{-1} \sum
_{l=1}^{2^{J_{n}}} \E\bigl[ \phi_{l}(W)^{2}
\bigl\{(g_{0} - P_{J_{n}}g_{0}) (X)\bigr
\}^{2} \bigr]
\\
&&\qquad= n^{-1} \sum_{l=1}^{2^{J_{n}}} \iint
\phi_{l}(w)^{2} \bigl\{(g_{0} -
P_{J_{n}}g_{0}) (x)\bigr\}^{2} f_{X,W}(x,w)
\,dx \,dw
\\
&&\qquad\leq n^{-1} C_{1} \| g_{0} - P_{J_{n}}g_{0}
\|^{2} \times\sum_{l=1}^{2^{J_{n}}} \int
\phi_{l}(w)^{2} \,dw =o\bigl(2^{J_{n}}/n\bigr).
\end{eqnarray*}
Hence
\begin{eqnarray*}
\En\bigl[ \phi^{J_{n}}(W_{i}) R_{i} \bigr] &=& \En
\bigl[ \phi^{J_{n}} (W_{i}) U_{i} \bigr]
\\
&&{} + \E\bigl[ \phi^{J_{n}}(W) (g_{0}-P_{n}g_{0})
(X)\bigr] + \mathrm{Rem}
\end{eqnarray*}
with $\| \mathrm{Rem} \|_{\ell^{2}} = o_{P}(\sqrt{2^{J_{n}}/n})$. The
second term on the right-hand side is $O(\tau_{J_{n}}2^{-J_{n}s})$ in
the Euclidean norm.
Together with Lemma \ref{lemprf3}, we have
\[
\Prob\bigl\{ \bigl\| \En\bigl[ \phi^{J_{n}}(W_{i})
R_{i} \bigr] \bigr\|^{2}_{\ell^{2}} > D_{1} \bigl(
\tau_{J_{n}}^{2} 2^{-2J_{n}s} + 2^{J_{n}}/n\bigr)
\bigr\} \to0.
\]

Moreover, by Lemma \ref{lemB1}, we have
\[
\tr(\hat{\Phi}_{WW}) \leq2^{J_{n}} s_{\max}(\hat{
\Phi}_{WW}) \leq C_{1}\bigl(1+o_{P}(1)
\bigr)2^{J_{n}}.
\]
Taking these together, we have
\[
\Prob\bigl\{ \bigl\| \En\bigl[ \phi^{J_{n}}(W_{i})
R_{i} \bigr] \bigr\|_{\ell^{2}}^{2} + n^{-1} \tr(
\hat{\Phi}_{WW}) \leq D_{2}\bigl( \tau_{J_{n}}^{2}
2^{-2J_{n}s} + 2^{J_{n}}/n\bigr) \bigr\} \to1.
\]
By the proof of Theorem \ref{thm1}, this leads to the desired conclusion.
\end{pf*}

\section{Discussion}\label{sec6}

We have studied the asymptotic properties of quasi-poste\-rior
distributions against sieve priors in the NPIV model and given some
specific priors for which the quasi-posterior distribution (the
quasi-Bayes estimator) attains the minimax optimal rate of contraction
(convergence, resp.). These results greatly sharpen the previous
work \cite{LJ11}. We end this paper with two additional discussions.

\subsection{Multivariate case}\label{sec6.1}
In this paper, we have focused on the case where $X$ and $W$ are
scalar, mainly to avoid the notational complication. It is not
difficult to see that the results naturally extend to the case where
$X$ and $W$ are vectors with the same dimension, by considering tensor
product sieves (the contraction/convergence rates will then deteriorate
as the dimension grows). We can also consider the following more
general situation as in Section 3 of \cite{HH05}: suppose that $Y$ is a
scalar random variable, $X$ and $W$ are random vectors with the same
dimension, and $Z$ is another random vector (whose dimension may be
different from $X$), and suppose that we are interested in estimating
the function $g_{0}$ identified by the conditional moment restriction:
$\E[ Y \mid Z,W ] = \E[ g_{0} (X,Z) \mid Z, W]$ or $Y=g_{0}(X,Z) + U$
with $\E[ U \mid Z,W ] = 0$ (i.e., $X$ and $Z$ are endogenous and
exogenous explanatory variables, resp.). In principle, the
analysis can be reduced to the case where there are no exogenous
variables by conditioning on $Z=z$ (so the sieve measure of
ill-posedness can be defined by the one conditional on $Z=z$). More
precisely, when $Z$ is discretely distributed with finitely many mass
points, then $g_{0}(x, z)$, where $z$ is a mass point, can be estimated
by using only observations $i$ for which $Z_{i}=z$. When $Z$ is
continuously distributed, then $g_{0}(x,z)$ can be estimated by using
observations $i$ for which $Z_{i}$ is ``close'' to $z$; one way is to
use kernel weights as in Section 4.2 of \cite{H11b}. However, the
detailed analysis of this case is not presented here for brevity.

\subsection{Direction of future research}\label{sec6.2}
Finally, we make some remarks on the direction of future research. First,
as also noted by \cite{LJ11}, (adaptive) selection of the resolution
level $J_{n}$ in a (quasi-)Bayesian or ``empirical'' Bayesian approach
is an important topic to be investigated. Second, a (quasi-)Bayesian
analysis is typically useful in the analysis of complex models in
which frequentist estimation is difficult to implement due to
nondifferentiability/nonconvex nature of loss functions. This
usefulness comes from the fact that a (quasi-)Bayesian approach is
typically able to avoid numerical optimization. See \cite{CH03} and
\cite{LTW07} for the finite-dimensional case. In infinite-dimensional
models, such a computational challenge in frequentist estimation occurs
in the analysis of nonparametric instrumental quantile regression
models \cite{HL07,CP11,GS11}. In that model, a typical loss function
contains the indicator function and hence highly nonconvex. In such a
case, the computation of an optimal solution is by itself difficult,
and a solution obtained, if possible, is typically not guaranteed to be
globally optimal since there may be many local optima. It is hence of
interest to extend the results of the paper to nonparametric
instrumental quantile regression models. The extension to the quantile
regression case, which is currently under investigation, is highly
nontrivial since the problem of estimating the structural function
becomes a \emph{nonlinear} ill-posed inverse problem and a delicate
care of the stochastic expansion of the criterion function is needed.

\section*{Acknowledgments}

A major part of the work was done while the author was visiting the
Department of Economics, MIT. He would like to thank Professor Victor
Chernozhukov for his suggestions and encouragements, as well as
Professor Yukitoshi Matsushita for his constructive comments. Also he
would like to thank the Editor, Professor Runze Li, the Associate Editor, and
anonymous referees for their insightful comments that helped improve
on the quality of the paper.

\begin{supplement}
\stitle{Supplement to ``Quasi-Bayesian analysis of
nonparametric instrumental variables models''}
\slink[doi]{10.1214/13-AOS1150SUPP} 
\sdatatype{.pdf}
\sfilename{aos1150\_supp.pdf}
\sdescription{This supplemental file contains the additional technical
proofs omitted in the main text, and some technical tools used in the
proofs.}
\end{supplement}


\printaddresses

\end{document}